\RequirePackage{algorithmicx}
\RequirePackage{algpseudocode}
\documentclass[final]{siamart220329}  

\usepackage{siampaper}
\usepackage{defs}

\newcommand{\cahiernumber}{65}  

\def\thistitle{Nonsmooth exact penalty methods for equality-constrained optimization:\\ complexity and implementation}
\def\authorone{Youssef Diouane}
\def\authortwo{Maxence Gollier}
\def\authorthree{Dominique Orban}

\hypersetup{
  pdftitle={\thistitle},
  pdfauthor={\authorone{}, \authortwo{} and \authorthree{}},
  pdfsubject={Report Subject},
  pdfkeywords={equality-constrained optimization, penalty methods, proximal methods},
}

\title{\thistitle}
\author{%
  \authorone\thanks{%
    GERAD and Department of Mathematics and Industrial Engineering, Polytechnique Montr\'eal.
    E-mail: \mailto{youssef.diouane@polymtl.ca}.
    Research supported by an NSERC Discovery grant.
 }
  \and
  \authortwo\thanks{%
    GERAD and Department of Mathematics and Industrial Engineering, Polytechnique Montr\'eal.
    E-mail: \mailto{maxence-2.gollier@polymtl.ca}.
 }
  \and
  \authorthree\thanks{%
    GERAD and Department of Mathematics and Industrial Engineering, Polytechnique Montr\'eal.
    E-mail: \mailto{dominique.orban@gerad.ca}.
    Research supported by an NSERC Discovery grant.
 }
}

\date{\today}

\begin{document}

\maketitle

\thispagestyle{firstpage}
\pagestyle{myheadings}

\begin{abstract}
  Penalty methods are a well known class of algorithms for constrained optimization.
  They transform a constrained problem into a sequence of unconstrained \emph{penalized} problems in the hope that approximate solutions of the latter converge to a solution of the former.
  If Lagrange multipliers exist, exact penalty methods ensure that the penalty parameter only need increase a finite number of times, but are typically scorned in smooth optimization for the penalized problems are not smooth.
  This led researchers to consider the implementation of exact penalty methods inconvenient.
  Recent advances in proximal methods have led to increasingly efficient solvers for nonsmooth optimization.
  We study a general exact penalty algorithm and use it to show that the exact $\ell_2$-penalty method for equality-constrained optimization can, in fact, be implemented efficiently by solving the penalized problem using a proximal-type algorithm.
  We study the convergence of our algorithm and establish a worst-case complexity bound of $\mathcal{O}(\epsilon^{-2})$ to bring a stationarity measure below \(\epsilon > 0\) under the Mangarasian-Fromowitz constraint qualification and Lipschitz continuity of the objective gradient and constraint Jacobian.
  While the Lipschitz continuity of the objective gradient is not required for convergence in view of recent works, it is used in our analysis to derive the complexity bound.
  In a degenerate scenario where the penalty parameter grows unbounded, the complexity becomes \(\mathcal{O}(\epsilon^{-8})\), which is worse than another bound found in the literature.
  We justify the difference by arguing that our feasibility measure is properly scaled\@.
  Finally, we report numerical experience on small-scale problems from a standard collection and compare our solver with an augmented-Lagrangian and an SQP method.
  Our preliminary implementation is superior to the augmented Lagrangian in terms of robustness and efficiency, and is competitive with the SQP method in terms of robustness, though the latter retains a slight edge in terms of number of problem function evaluations.
\end{abstract}



\section{Introduction}%
\label{sec:introduction}

We consider the problem
\begin{equation}%
  \label{eq:nlp}
  \minimize{x \in \R^n} \ f(x) \quad \st \ c(x) = 0,
\end{equation}
where \(f: \R^n \to \R\) and \(c: \R^n \to \R^m\) are \(\mathcal{C}^1\) and both may be nonconvex.
We solve~\eqref{eq:nlp} by solving a sequence of unconstrained, nonsmooth penalized problems
\begin{equation}%
  \label{eq:penalty-nlp}
  \minimize{x \in \R^n} \ f(x) + \tau h(c(x)),
\end{equation}
where \( h = \|\cdot\| \) denotes a norm on \( \R^n \), and \( \tau > 0 \) is the penalty parameter.
This approach was first proposed by \citet{pietrzykowski-1969}, who used the $\ell_1$-norm.
Penalty approaches such as~\eqref{eq:penalty-nlp} are attractive because, under standard assumptions, for all sufficiently large and finite $\tau$, solutions of~\eqref{eq:nlp} are solutions of~\eqref{eq:penalty-nlp}.
However, the nonsmoothness of \(h\) caused them to fall out of favor in practice, and other methods, such as the augmented-Lagrangian method \citep{hestenes-1969,powell-1969} were preferred.

Independently, attention was given recently to nonsmooth regularized problems
\begin{equation}%
  \label{eq:nonsmooth-nlp}
  \minimize{x \in \R^n} \ f(x) + h(x),
\end{equation}
where $f$ is as before, and $h$ is proper and lower semi-continuous (lsc), and may be nonconvex.
In particular, proximal methods \citep{lions-mercier-1979,fukushima-mine-1981} have been studied intensively in the last decade and have led to increasingly efficient implementations.

We show that proximal methods can be used to implement exact penalty methods so they perform similarly to augmented-Lagrangian approaches.
To do so, we study an intuitive algorithm that enacts the exact penalty scheme, but in which subproblems are solved by a first-order proximal method proposed in \citep{aravkin-baraldi-orban-2022,aravkin-baraldi-leconte-orban-2021}, and a second-order proximal method proposed in \citep{diouane-habiboullah-orban-2024}.
We establish convergence of an appropriate stationarity measure to zero, and a worst-case complexity bound of $\mathcal{O}(\epsilon^{-2})$ to bring said measure below \(\epsilon > 0\) under the Mangarasian-Fromowitz constraint qualification and Lipschitz continuity of the objective gradient and constraint Jacobian.
While the Lipschitz continuity of the objective gradient is not required for convergence in view of \citep{diouane-habiboullah-orban-2024}, it is used in our analysis to derive the complexity bound.
In a degenerate scenario where the penalty parameter grows unbounded, the complexity becomes \(\mathcal{O}(\epsilon^{-8})\), which is worse than the bound \(\mathcal{O}(\epsilon^{-5})\) developed in \citep[\S\(3\)]{cartis-gould-toint-2011}.
Although both bounds are correct, we justify the difference by arguing that our feasibility measure is properly scaled.
We use our general analysis to implement the $\ell_2$-norm penalty method, though $\ell_1$ and $\ell_\infty$ variants may also be implemented efficiently, building upon the same framework.
As far as we are aware, ours is the first implementation of the exact penalty scheme based on proximal methods.
We report numerical experience on small-scale problems from a standard collection and compare our solver with an augmented-Lagrangian and an SQP method.
Our preliminary implementation is superior to the augmented Lagrangian in terms of robustness and efficiency, and is competitive with the SQP method in terms of robustness, though the latter retains a slight edge in terms of number of problem function evaluations.

The paper is organized as follows.
In \Cref{sec:background} we recall key background concepts and results.
In \Cref{sec:models} we introduce the models of \(f\) and \(c\) in~\eqref{eq:nlp} that are used to compute steps in our algorithm.
In \Cref{sec:algorithm-convergence}, we present our algorithm and derive complexity bounds.
In \Cref{sec:prox-operators} we derive procedures to evaluate the proximal operators that arise in our algorithm.
Additionally, we provide algorithms to evaluate the required proximal operators efficiently in practice.
In \Cref{sec:numerical} we show numerical results and experiments conducted on a variety of problems.
Finally, we provide a closing discussion in \Cref{sec:discussion}.

\subsection*{Related research}

\citet{pietrzykowski-1969} first reported the advantage of using a nonsmooth penalty function.
His analysis only covered the $\ell_1$-penalty but subsequently, \citet{charalambous-1978}, \citet{han-mangasarian-1979}, \citet{coleman-conn-1980}, \citet{bazaraa-goode-1982}, and \citet{huang-ng-1994} strengthened his results.

\citet{cartis-gould-toint-2011} analyze general \emph{composite} objectives, i.e., of the form $h(c(x))$, where $h$ is convex and globally Lipschitz using a trust-region approach and a quadratic regularization variant.
They then apply their results to an exact penalty method similar to ours.
\citet{grapiglia-yuan-yuan-2015} generalize the results of \citep{cartis-gould-toint-2011} to various optimization frameworks under similar assumptions.
Contrary to their algorithms, ours is concrete in the sense that we provide precise guidance on how to compute steps.
Indeed, we show that solving the subproblem of the exact penalty method can be done by solving a trust-region subproblem for which algorithms are well-known for certain choices of norms.
By contrast, the algorithms in \citep{grapiglia-yuan-yuan-2015,cartis-gould-toint-2011} require to solve a nonsmooth subproblem with a trust-region constraint, which is difficult in practice, and those authors provide no indication as to how that problem may be solved.
Additionally, our convergence analysis holds under weaker assumptions.

Several authors reformulate~\eqref{eq:penalty-nlp} as a smooth inequality-constrained problem.
For instance,~\eqref{eq:penalty-nlp} can be written equivalently
\[
  \minimize{x \in \R^n, u \in \R^m} \ f(x) + \tau u^T e \quad \st \ -u \leq c(x) \leq u,
\]
when \(h(x) = \|x\|_1\), where \(e\) is the vector of ones.
That is what \citet{gould-orban-toint-2015b} do before applying an interior-point method.
The literature review of \citep{gould-orban-toint-2015b} provides a number of related references.
A similar transformation exists for \(h(x) = \|x\|_\infty\), or any other polyhedral norm.

\citet{estrin-friedlander-orban-saunders-2020a,estrin-friedlander-orban-saunders-2020b} implement a smooth exact penalty function originally proposed by \citet{fletcher-1970}.
They challenge the notion that evaluating the penalty function and its gradient is costlier than in other widely-accepted constrained frameworks, and detail an efficient implementation.

Originally proposed by \citet{hestenes-1969} and \citet{powell-1969}, the augmented-Lagrangian approach may be defined as a quadratic penalty applied to the Lagrangian, and was viewed by Powell as a shifted quadratic penalty function.
Two of its attractive features are that it is smooth, and it acts as an exact penalty once optimal multipliers have been identified.
\citet{bertsekas-1982} analyzes its convergence thoroughly and improves the original results.
We refer interested readers to \citep{bertsekas-1982} for further information.
Renowned implementations of augmented-Lagrangian methods include LANCELOT \citep{conn-gould-toint-1991,conn-gould-toint-1992}, MINOS \citep{murtagh-saunders-1978,murtagh-saunders-1982}, and ALGENCAN \citep{andreani-birgin-martinez-schuverdt-2008,birgin-martinez-2014}.
In our numerical experiments, we use the recent implementation of \citet{percival-jl}, named Percival.
In the context of exact penalty methods, \citet{kuhlmann-buskens-2018} proposed an augmented-Lagrangian-like approach in which the quadratic penalty is replaced by an exact $\ell_2$-penalty.
Our approach differs in that we add the penalty term directly to the objective of~\eqref{eq:nlp}.
Moreover, while both approaches lead---via the exact penalty term in~\eqref{eq:penalty-nlp}---to optimality conditions involving a trust-region-like term, they ignore this term when computing steps, using it only to update penalty parameters, which is in contrast with our approach.
Additionally, their approach consists in applying Newton's method to the first-order optimality condition of the penalized problem, which differs from the proximal-gradient method that we use.
Their algorithm has been implemented within the solver WORHP \citep{buskens-wassel-2013}.

The proximal-gradient (PG) method \citep{fukushima-mine-1981,lions-mercier-1979} aims to solve problems of the form~\eqref{eq:nonsmooth-nlp} where $f$ is \(\mathcal{C}^1\) and $h$ is proper and lsc.
In that context, \(h\) often has regularizing power; it promotes solutions with desirable features, such as sparsity.
Variants in the literature mainly differ in the assumptions on \(f\) and \(h\).
Numerous authors restrict their attention to convex \(f\) and/or \(h\).
\citet{parikh-boyd-2013} review PG in the convex case with an insightful chapter on their interpretation.
\citet{lee-sun-saunders-2014} study the convergence of a proximal Newton method in which steps are computed with PG where both \(f\) and \(h\) are convex.
\citet{bolte-sabach-teboulle-2013} present a method for objectives of the form $g(x) + Q(x,y) + h(y)$ where both $g$ and $h$ are proper, lsc, and $Q$ is \(\mathcal{C}^1\).
\citet{tseng-2000} proposes a general framework for accelerated PG\@.
Proximal algorithms are related to augmented-Lagrangian methods \citep{rockafellar-1973,rockafellar-1976}.
To the best of our knowledge, almost all proximal methods in the literature require evaluating a proximal operator---which involves the solution of a nonsmooth problem consisting of the sum of a squared norm and \(h\)---at each iteration.
Evaluating that operator for \(h \circ c\) in~\eqref{eq:penalty-nlp} would be impractical and prohibitively expensive.
To compute a step, we use a variant of PG with adaptive step size developed by \citet{aravkin-baraldi-orban-2022} that does not have convexity restrictions.
Crucially, their method allows for a model of \(h\) to be used at each iteration, a feature that makes evaluating the proximal operator feasible.
In a second stage, we use the recent proximal quasi-Newton method of \citet{diouane-habiboullah-orban-2024}, which may be seen as a natural generalization of the method of \citep{aravkin-baraldi-orban-2022}.

\subsection*{Notation}

The identity matrix of size \(n\) is \(I_n\), or \(I\) if the context is clear.
We use \(J(x)\) to denote the Jacobian of \(c\) at \(x\).
We say a function \(f\) is \(\mathcal{C}^1\) if it is continuously differentiable.
For any matrix \(A\), its Moore-Penrose pseudo-inverse is \(A^\dagger\).
For \(p \geq 1\), the vector $\ell_p$-norm is \(\| \cdot \|_p\).
\(\| \cdot \|\) is any norm on \(\R^n\) and \(\|\cdot\|_*\) is its dual norm, defined as \(\|y\|_* \coloneqq \sup_{\|x\| \leq 1} x^T y\).
For matrices, \(\| \cdot \|_p\) represents the operator norm with respect to the vector \(\ell_p\)-norm.
Similarly, $\B_p$ represents the unit ball, centered at the origin in the $\ell_p$-norm.
For any $\Delta > 0$, $\Delta\B_p$ represents the ball of radius $\Delta > 0$ centered at the origin.
If a set $\mathcal{S}$ is finite, $|\mathcal{S}|$ represents its cardinality.
For any set $\mathcal{S} \subseteq \R^n$, $\chi(\cdot \mid \mathcal{S})$ is the \textit{indicator} function, namely, for any $x \in \R^n$, $\chi(x \mid \mathcal{S}) = 0$ if $x \in \mathcal{S}$ and $+ \infty$ otherwise.
If $g$, $h$ are two positive functions of $\epsilon > 0$, the notation $g = \mathcal{O}(h)$ means that there is \(C > 0\) such that \(\limsup_{\epsilon \to 0} g(\epsilon)/h(\epsilon) \leq C\).
The notation \(g = o(h)\) means that \(\lim_{\epsilon \to 0} g(\epsilon)/h(\epsilon) = 0\).

\section{Background}%
\label{sec:background}

The Mangasarian-Fromowitz constraint qualification (MFCQ) holds at \(x \in \R^n\) for~\eqref{eq:nlp} if $J(x)$ has full row rank.

The element $\bar{x} \in \R^n$ is a strict minimum of~\eqref{eq:nlp} if \(c(\bar{x}) = 0\) and there is an open set $\mathcal{V}$ containing $\bar{x}$ such that \(f(\bar{x}) < f(x)\) for all \(x \in \mathcal{V}\) satisfying \(c(x) = 0\).

Exactness of the norm penalty relies on existence of Lagrange multipliers, and means the following.

\begin{proposition}[\protect{\citealp[Theorem 4.4]{han-mangasarian-1979}}]%
  \label{th:exact-penanlty}
  If \(\bar{x}\) is a strict minimum of~\eqref{eq:nlp} where the MFCQ holds, for every $\tau \geq \|\bar{y}\|_*$, $\bar{x}$ is a local minimum of~\eqref{eq:penalty-nlp}, where $\bar{y}$ is the unique vector of Lagrange multipliers at $\bar{x}$.
\end{proposition}

Consider $\phi: \R^n \rightarrow \widebar{\R}$ and $\bar{x} \in \R^n$ where $\phi(\bar{x}) \in \R$.
We say that $v \in \R^n$ is a regular subgradient of $\phi$ at $\bar{x}$, and we write $v \in \hat{\partial} \phi(\bar{x})$, if \(\phi(x) \geq \phi(\bar{x}) + v^T (x - \bar{x}) + o(\|x - \bar{x}\|_2)\).
The set of regular subgradients is called the Fréchet subdifferential.
We say that $v$ is a general subgradient of $\phi$ at $\bar{x}$, and we write $v \in \partial\phi (\bar{x})$, if there are $\{x_k\}$ and $\{v_k\}$ such that $\{x_k\} \to \bar{x}$, $\{\phi(x_k)\} \to \phi(\bar{x})$, $v_k\ \in \hat{\partial}\phi(x_k)$ for all \(k\) and $\{v_k\} \to v$.
The set of general subgradients is called the limiting subdifferential \citep[Definition 8.3]{rockafellar-wets-1998}.

We call  $h: \R^n \rightarrow \widebar{\R}$ proper if $h(x) > -\infty$ for all $x \in \R^n$ and $h(x) < +\infty$ for at least one $x$.
\(h\) is lsc at $\bar{x}$ if $\liminf_{x \to \bar{x}} h(x) = h(\bar{x})$.

If $\phi$ is \(\mathcal{C}^1\) over a neighborhood of $x$, $\partial\phi(x) = \{\nabla \phi(x)\}$ \citep[\S\(8.8\)]{rockafellar-wets-1998}.
In what follows, we rely on the following criticality property.

\begin{proposition}[\protect{\citealp[Theorem 10.1]{rockafellar-wets-1998}}]{}%
  \label{prop:subdifferential}
  If $\phi: \R^n \rightarrow \R$ is proper and has a local minimum at $\bar{x}$, then $0 \in \hat{\partial} \phi(\bar{x}) \subseteq \partial \phi(\bar{x})$.
  If $\phi = f + h$, where f is \(\mathcal{C}^1\) over a neighborhood of $\bar{x}$ and h is finite at $\bar{x}$, then $\partial \phi(\bar{x}) = \nabla f(\bar{x}) + \partial h(\bar{x})$.
\end{proposition}

The element $\bar{x} \in \R^n$ is a KKT point of~\eqref{eq:nlp} if there is $\bar{y} \in \R^m$ such that \(\nabla f(\bar{x}) = J{(\bar{x})}^T \bar{y}\) and \(c(\bar{x}) = 0\).
Moreover, if there exists \(\epsilon_1, \epsilon_2 > 0\) such that there is a norm for which \(\|\nabla f(\bar{x}) - J{(\bar{x})}^T \bar{y}\| \leq \epsilon_1\) and \(\|c(\bar{x})\| \leq \epsilon_2\), then we call \(\bar{x}\) an \((\epsilon_1, \epsilon_2)\)-approximate KKT point of~\eqref{eq:nlp}.
Since all norms are equivalent on \(\R^n\), this definition of an approximate KKT point is independent of the specific norm chosen.

For proper lsc $h$ and step length $\nu > 0$, the proximal mapping of $\nu h$ at \(x\) is
\begin{equation}%
  \label{def:moreau-prox}
  \prox{\nu h}(x) \coloneqq \argmin{u} \tfrac{1}{2} \nu^{-1} \|u-x\|^2_2 + h(u).
\end{equation}

\begin{proposition}[\protect{\citealp[Theorem 1.25]{rockafellar-wets-1998}}]%
  \label{prop::Moreau}
  Let $h: \R^n \rightarrow \widebar{\R}$ be proper lsc and bounded below.
  For all $\nu > 0$ and all $x \in \R^n$, $\prox{\nu h}(x)$ is nonempty and compact.
\end{proposition}

The proximal mapping of \(\nu h\) is nonempty under more general conditions than stated in \Cref{prop::Moreau}, e.g., prox-boundedness \citep[Definition~\(1.23\)]{rockafellar-wets-1998}, but boundedness is sufficient for our purposes.

We now intentionally use the notation \(\varphi\) and \(\psi\) instead of \(f\) and \(h\) for reasons that become clear in subsequent sections.
Inspired by~\eqref{eq:penalty-nlp}, we consider the generic nonsmooth problem
\begin{equation}%
  \label{eq:proximal-gradient-method}
  \minimize{s} \varphi(s) + \psi(s),
\end{equation}
where $\varphi$ is \(\mathcal{C}^1\), and $\psi$ is proper, lsc, and bounded below.
A popular method to solve~\eqref{eq:proximal-gradient-method} is the \textit{proximal-gradient} (PG) method \citep{fukushima-mine-1981,lions-mercier-1979}.
The PG iteration is
\begin{equation}%
  \label{eq:prox-iterate}
  s_{i+1} \in \underset{\nu_i \psi}{\text{prox}}(s_i - \nu_i \nabla\varphi(s_i)), \quad i \geq 0,
\end{equation}
where $\nu_i > 0$ is a step length.
Descent is guaranteed when \(\nabla \varphi\) is Lipschitz-continuous and \(\nu_i\) is chosen appropriately.
Recall that $\phi: \R^{n_1} \rightarrow \R^{n_2}$ is Lipschitz-continuous with Lipschitz constant $L > 0$ whenever \(\|\phi(x) - \phi(y)\|_2 \leq L\|x - y\|_2\) for all \(x\), \(y \in \R^{n_1}\).

\begin{proposition}[\protect{\citealp[Lemma 2]{bolte-sabach-teboulle-2013}}]%
  \label{prop:prox-decrease}
  Let $\varphi$ be \(\mathcal{C}^1\), $\nabla\varphi$ be Lipschitz-continuous with Lipschitz constant $L \geq 0$, and $\psi$ be proper, lsc and bounded below.
  For any $0 < \nu < L^{-1}$ and $s_0 \in \R^n$ where $\psi$ is finite,~\eqref{eq:prox-iterate} is such that \((\varphi + \psi)(s_{i+1}) \leq (\varphi + \psi)(s_i) - \tfrac{1}{2} (\nu^{-1} - L) \| s_{i+1} - s_i \|^2_2\) for all \(i \geq 0\).
  If \(L = 0\), such as arises when \(\varphi\) is linear, the above is taken to mean \(L^{-1} = +\infty\).
\end{proposition}

\section{Models}%
\label{sec:models}

For given $x \in \R^n$, consider models
\begin{subequations}%
  \label{eq:models}
  \begin{align}%
    \varphi(s;x) & \approx f(x + s),                     \\
    \psi(s;x)    & \approx h(c(x + s)) = \| c(x + s) \|, \\
    m(s;x)       & \coloneqq \varphi(s;x) + \psi(s;x),
  \end{align}
\end{subequations}
where \(\approx\) means that the left-hand side is an approximation of the right-hand side.
The type of approximation expected of~\eqref{eq:models} is formalized in~\cref{modelassumption:1}.
\begin{modelassumption}%
  \label{modelassumption:1}
  For any $x \in \R^n$, $\varphi(\cdot;x)$ is \(\mathcal{C}^1\), and satisfies $\varphi(0;x) = f(x)$ and $\nabla_s \varphi(0;x) = \nabla f(x)$.
  For any $x \in \R^n$, $\psi$ is proper lsc, and satisfies $\psi(0;x) = h(c(x))$ and $\partial_s \psi(0;x) = \partial (h \circ c)(x)$.
\end{modelassumption}
In \Cref{modelassumption:1}, we use the notation \(\partial_s\) to emphasize that the subdifferential of \(\psi(s; x)\) is computed with respect to its variable \(s\) while keeping \(x\) fixed.

We first focus our attention on the first order models
\begin{subequations}%
  \label{eq:modLinear}
  \begin{align}%
    \varphi(s; x) & \coloneqq f(x) + \nabla {f(x)}^T s,
    \label{eq:modphi1}%
    \\
    \psi(s; x)    & \coloneqq h(c(x) + J(x)s) = \| c(x) + J(x)s \|.
    \label{eq:modpsi}
  \end{align}
\end{subequations}
For future reference, define
\begin{equation}%
  \label{eq:g}
  g_x  : \R^n \to \R^m, \quad g_x(s)  \coloneq c(x) + J(x) s \quad (x \in \R^n),
\end{equation}
From \citep[Section D.3]{hiriart-urruty-lemaréchal-2001},
\begin{equation}%
  \label{eq:norm-subdif}
  \partial h(z) = \{ \|y\|_* \leq 1 \mid  y^T z = \|z\|\}.
\end{equation}
It follows from~\eqref{eq:norm-subdif} that for any $y \in \partial h(z)$,
\begin{equation}
  \label{eq:norm-subdif-zero}
  \|y\|_* < 1 \Longrightarrow z = 0,
\end{equation}
because, from H\"older's inequality and~\eqref{eq:norm-subdif}, we have $\|z\| = y^T z \leq \|y\|_* \|z\|$
The following lemma gives the subdifferential of $\psi$, which will be useful to derive a stationarity measure for~\eqref{eq:penalty-nlp}.

\begin{lemma}%
  \label{lemma:subdif}
  Let $x \in \R^n$.
  Then, for all \(s \in \R^n\),
  \begin{equation}%
    \label{eq:psi-subdif}
    \partial_s \psi(s; x) = {J(x)}^T \partial h(c(x) + J(x)s)
  \end{equation}
  where \(\psi\) is defined in~\eqref{eq:modpsi} and $\partial h$ is given in~\eqref{eq:norm-subdif}.
\end{lemma}

\begin{proof}%
  Let $g_x$ be as in~\eqref{eq:g}.
  Since $\dom h = \R^m$ and $c(x) \in \Range(g_x) \subseteq \R^m$, the result follows from \citep[Theorem 23.9]{rockafellar-1970} and~\eqref{eq:modpsi}.
\end{proof}

Since \(h\) is bounded below, lsc and convex, \citep[Theorem 10.6]{rockafellar-wets-1998} yields \(\partial (h \circ c)(x) = {J(x)}^T \partial h(c(x))\), which, by \Cref{lemma:subdif}, is equal to \(\partial_s \psi(0; x)\).
Thus,~\eqref{eq:modLinear} satisfy \Cref{modelassumption:1}.

\section{Algorithm and convergence analysis}%
\label{sec:algorithm-convergence}

At iteration $k \in \N$, we solve~\eqref{eq:penalty-nlp} inexactly for fixed $\tau_k > 0$.
The process of updating \(\tau_k\) is the \textit{outer} iteration.
Solving~\eqref{eq:penalty-nlp} with \(\tau = \tau_k\) is the \(k\)-th set of \textit{inner} iterations with inner iterates $x_{k,j}$.

We solve~\eqref{eq:penalty-nlp} with method R2 \citep[Algorithm~\(6.1\)]{aravkin-baraldi-orban-2022}; a quadratic regularization method that uses models~\eqref{eq:modLinear} and may be viewed as an adaptive variant of PG that does not require knowledge of the Lipschitz constant of the gradient of \(f\) (or even its existence for convergence).
At inner iteration $j$, we compute
\begin{subequations}%
  \label{eq:def-skj1}
  \begin{align}
    \{ s_{k, j, \mathrm{cp}} \} & = \argmin{s} \varphi(s ; x_{k, j}) + \tau_k\psi(s ; x_{k, j}) + \tfrac{1}{2} \sigma_{k, j} \|s\|^2_2
    \label{eq:prox-problem}
    \\
                                & = \prox{\sigma_{k,j}^{-1} \tau_k\psi(\cdot ; x_{k,j})}(-\sigma_{k,j}^{-1} \nabla \varphi(0 ; x_{k,j})),
    \label{eq:proximal-step-1}
  \end{align}
\end{subequations}
where $\sigma_{k,j} > 0$ is a regularization parameter, i.e., $s_{k, j, \mathrm{cp}}$ is the first step of PG for the minimization of $\varphi(s ; x_{k,j}) + \tau_k\psi(s ; x_{k,j})$ with step length $\sigma_{k,j}^{-1}$ initialized with $s_{k,j,0} = 0$ \citep[\S\(3.2\)]{aravkin-baraldi-orban-2022}.
Convexity of $\varphi(\cdot; x_{k, j})$ and $\tau_k\psi(\cdot; x_{k, j})$ ensures that the proximal operator in~\eqref{eq:def-skj1} is single-valued.
The subscript ``cp'' stands for \emph{Cauchy point}, of which \(s_{k, j, \mathrm{cp}}\) is an appropriate generalization to the nonsmooth context.
We then set either $x_{k,j+1} \coloneqq x_{k,j} + s_{k, j, \mathrm{cp}}$ or $x_{k,j+1} \coloneqq x_{k,j}$ depending on whether \(s_{k, j, \mathrm{cp}}\) results in sufficient decrease in the objective of~\eqref{eq:penalty-nlp} or not.
The complete algorithm is given in \Cref{sec:inner_algs}; we refer the reader to \citep{aravkin-baraldi-orban-2022,aravkin-baraldi-leconte-orban-2021} for details.
In particular, the quantity
\begin{equation}%
  \label{eq:xi}
  \xi(x_{k,j} ; \sigma_{k,j}, \tau_k) \coloneqq f(x_{k,j}) + \tau_k h(c(x_{k,j})) - (\varphi + \tau_k\psi)(s_{k, j, \mathrm{cp}} ; x_{k,j})
\end{equation}
is key to defining the stationarity measure \citep{aravkin-baraldi-orban-2022,aravkin-baraldi-leconte-orban-2021}
\begin{equation}%
  \label{eq:xi-half}
  \sigma_{k,j}^{1/2} \xi{(x_{k,j} ; \sigma_{k,j},\tau_k)}^{1/2}.
\end{equation}

\begin{proposition}%
  \label{prop:kkt}
  Let $\sigma > 0$, $\tau > 0$, \(\kappa_c > 0\), \(\epsilon > 0\) and \(\bar{x} \in \R^n\).
  Assume that \(\sigma^{1/2}\xi{(\bar{x}; \sigma, \tau)}^{1/2} \leq \epsilon\), and let \(\{s_{\mathrm{cp}}\} = \prox{\sigma^{-1}\tau\psi(\cdot; \bar{x})}(-\sigma^{-1}\nabla \varphi(0;\bar{x}))\).
  Then, there is $\bar{y} \in \tau \partial h(c(x) + J(x)s_{\mathrm{cp}})$ such that \(\| \nabla f(\bar{x}) + J{(\bar{x})}^T \bar{y} \|_2 \leq \sqrt{2}\epsilon\).
  Moreover, if $\|c(\bar{x})\|_2 \leq \kappa_c \epsilon$, \(\bar{x}\) is an \((\sqrt{2}\epsilon, \kappa_c \epsilon)\)-approximate KKT point of~\eqref{eq:nlp}.
\end{proposition}

\begin{proof}%
  By definition of \(s_{\mathrm{cp}}\) and~\eqref{eq:def-skj1},
  \begin{equation}%
    \label{eq:kkt1}
    \{ s_{\mathrm{cp}} \} = \argmin{s} \varphi(s; \bar{x}) + \tau \psi(s; \bar{x}) + \tfrac{1}{2} \sigma \| s \|^2_2.
  \end{equation}
  On the other hand, \Cref{prop:prox-decrease} and~\eqref{eq:xi} imply
  \begin{equation}%
    \label{eq:kkt2}
    \epsilon > \sigma^{1/2}\xi{(\bar{x}, \sigma, \tau)}^{1/2} \geq \frac{\sigma \|s_{\mathrm{cp}}\|_2}{\sqrt{2}}.
  \end{equation}
  Convexity, the first-order conditions of~\eqref{eq:kkt1} and \Cref{prop:subdifferential} yield
  \begin{equation}%
    \label{eq:kkt3}
    0 \in \nabla \varphi(s_{\mathrm{cp}}; \bar{x}) + \tau \partial_s \psi(s_{\mathrm{cp}}; \bar{x}) + \sigma s_{\mathrm{cp}}.
  \end{equation}
  Thus, by~\eqref{eq:modphi1}, \Cref{lemma:subdif} and injecting~\eqref{eq:kkt3} into~\eqref{eq:kkt2}, there is $y \in \partial h(c(x) + J(x)s_{\mathrm{cp}})$ such that \(\|\nabla f(\bar{x}) + \tau {J(\bar{x})}^T y\|_2 \leq \sqrt{2}\epsilon\), and the result holds with $\bar{y} \coloneq \tau y$.
  If additionally, $\|c(\bar{x})\|_2 \leq \kappa_c \epsilon$, \(\bar{x}\) is a \((\sqrt{2}\epsilon, \kappa_c \epsilon)\)-approximate KKT point of~\eqref{eq:nlp}.
\end{proof}

In view of \Cref{prop:kkt}, the aim of the outer iterations is to find a feasible point for~\eqref{eq:nlp}.
We focus on the feasibility problem
\begin{equation}%
  \label{eq:feasibility}
  \minimize{x \in \R^n} \ 0 \quad \st \ c(x) = 0.
\end{equation}
In the same fashion as we did for $\xi$, in the case where $f = 0$, we define a feasibility measure for the outer iterates based on
\begin{equation}%
  \label{eq:theta}
  \theta(x) \coloneq h(c(x)) - \psi(s^* ; x), \quad \{s^*\} = \prox{1 \psi(\cdot; x)} (0) = \argmin{s} \psi(s; x) + \tfrac{1}{2} \|s\|_2^2.
\end{equation}

Note that $\theta(x)$ is $\xi(x; 1, 1)$ in the case where $f=0$.
In particular, \Cref{prop:kkt} applies, so that for any \(\epsilon > 0\), \(\theta{(x)}^{1/2} \leq \epsilon\) implies \(\|{J(x)}^T y\|_2 \leq \sqrt{2}\epsilon\) for some \(y \in \partial h(c(x) + J(x)s^*)\).
From~\eqref{eq:norm-subdif} \(\|y\|_* \leq 1\) and, in particular, if \(\|y\|_* < 1\),~\eqref{eq:norm-subdif-zero} implies \(c(x) + J(x)s^* = 0\).
Hence,~\eqref{eq:theta} implies \(h(c(x)) \leq \epsilon^2\).

By analogy with~\eqref{eq:xi-half}, we define our feasiblity measure as \(\theta{(x)}^{1/2}\).
\Cref{alg:exactpen} summarizes the outer iteration.

\begin{algorithm}[ht]%
  \caption{%
    \label{alg:exactpen}
    Exact penalty algorithm.
  }
  \begin{algorithmic}[1]%

    \State{} Choose \(x_0 \in \R^n\), $\beta_1 > 0$ and $\tau_{0} > 0$.%

    \State{} Choose initial and final tolerances $\epsilon_0 \geq \epsilon > 0$, $0 < \beta_2 < 1$, and $\beta_{3}$, $\beta_{4} > 0$.%

    \For{$k = 0, 1, \ldots$}

    \State%
    \label{eq:exactpen_inner_stationnarity}
    Compute an approximate solution of~\eqref{eq:penalty-nlp} with $\tau \coloneq \tau_k$ starting from $x_k$ with initial step size $\sigma_{k,0} = \max(\beta_3 \tau_k, \, \beta_4)$ and minimal regularization parameter $\sigma_{k,\min} = \beta_4$.
    Stop at the first iteration $j$ such that
    \begin{equation*}%
      \sigma_{k,j}^{1/2} \xi{(x_{k,j} ; \sigma_{k,j},\tau_k)}^{1/2} \leq \epsilon_k.
    \end{equation*}

    \State{} Set $x_{k + 1} \coloneqq x_{k,j}$.%

    \State%
    \label{eq:exactpen_epsilon_sequence}
    If $\theta^{1/2}(x_{k+1}) > \epsilon_k$, choose $\tau_{k+1} \geq \tau_k + \beta_1$ and set $\epsilon_{k+1} \coloneqq \epsilon_k$.
    Otherwise, set $\tau_{k+1} \coloneqq \tau_k$ and $\epsilon_{k+1} \coloneqq \beta_2\epsilon_k$.

    \EndFor%
  \end{algorithmic}
\end{algorithm}

\Cref{alg:exactpen} is similar to \citep[Algorithm 3.1]{cartis-gould-toint-2011}, except that we do not use ``steering''; a procedure that guarantees that the outer and inner measures decrease at the same rate.
Instead, we use an increasing accuracy strategy of the inner measure on \Cref{eq:exactpen_epsilon_sequence}.

The following assumption ensures convergence of the inner iterations.

\begin{stepassumption}[\protect{\citealp[Step Assumption 6.1]{aravkin-baraldi-orban-2022}}]%
  \label{step_assumption:inner}
  For each outer iteration $k$ of \Cref{alg:exactpen}, there is $\kappa_{m,k} > 0$ such that for each corresponding inner iteration $j$, $s_{k, j, \mathrm{cp}}$ generated according to~\eqref{eq:def-skj1} satisfies
  \begin{equation}%
    |f(x_{k,j} + s_{k, j, \mathrm{cp}}) + \tau_k h(c(x_{k,j} + s_{k, j, \mathrm{cp}})) - (\varphi + \tau_k\psi)(s_{k, j, \mathrm{cp}} ; x_{k,j})| \leq \kappa_{m,k} \|s_{k, j, \mathrm{cp}}\|^2_2.
  \end{equation}
\end{stepassumption}

We make the following additional assumption on the growth of \(\kappa_{m,k}\) along the outer iterations.

\begin{stepassumption}%
  \label{step_assumption:kappa-tauk}
  There exist $\kappa_{m_1} > 0$ and $\kappa_{m_2} > 0$ such that for each outer iteration $k$ of \Cref{alg:exactpen}, $\kappa_{m,k} \leq \kappa_{m_1} \tau_k + \kappa_{m_2}$, where \(\kappa_{m,k} > 0\) is defined in \Cref{step_assumption:inner}.
\end{stepassumption}
\Cref{step_assumption:kappa-tauk} is satisfied when $f$ and $c$ have Lipschitz gradient and Jacobian, respectively.
Indeed, if $\nabla f$ and $J$ have Lipschitz constants $L_g$ and $L_J$, respectively, then for any $x$, $s \in \R^n$ and $\tau > 0$,
\(
|f(x + s) + \tau h(c(x + s)) - \varphi(s ; x) - \tau\psi(s ; x)| \leq \tfrac{1}{2} L_g \|s\|^2_2 + \tau \tfrac{1}{2} L_J \|s\|^2_2
\),
where we used the triangle inequality, the descent lemma for smooth functions \citep[Lemma 1]{bolte-sabach-teboulle-2013} on $f$ and $c$, and the fact that \(\tau h\) is Lipschitz-continuous with constant \(\tau\).
Given the recent developments by~\citet{diouane-habiboullah-orban-2024}, Lipschitz continuity of \(\nabla f\) is not necessary to establish convergence of \Cref{alg:exactpen}.
However, it remains essential to derive complexity bounds, which is the focus of the present analysis.
Informally, \Cref{step_assumption:inner} bounds the relative error of the models at each outer iteration, while \Cref{step_assumption:kappa-tauk} states that the model quality is expected to deteriorate at most linearly with the penalty parameter.

For a tolerance $\epsilon > 0$, we wish to find the number of inner and outer iterations of \Cref{alg:exactpen} until
\begin{equation}%
  \label{eq:stop-crit}
  {\theta(x_k)}%
  ^{1/2}%
  \leq \epsilon.
\end{equation}

We use the following notation in the complexity analysis:
\begin{align}%
  j_k(\epsilon_k) & \coloneq \min \{j \in \N \mid \sigma_{k,j}^{1/2} \xi{(x_{k,j} ; \sigma_{k,j},\tau_k)}^{1/2} \leq \epsilon_k \} \quad (k \in \N), \\
  k(\epsilon)     & \coloneq \min \{k \in \N \mid {\theta(x_k)}^{1/2} \leq \epsilon \}.
\end{align}
The total number of iterations we are looking for is then $\sum_{k = 0}^{k(\epsilon)} j_k(\epsilon_k)$.

The next result states \citep[Theorem~\(6.4\)]{aravkin-baraldi-orban-2022} with the bound given by \citet{aravkin-baraldi-leconte-orban-2021}.

\begin{theorem}[\protect{\citealp[Theorem~\(6.4\)]{aravkin-baraldi-orban-2022} and \citealp{aravkin-baraldi-leconte-orban-2021}}]%
  \label{th:inner-complexity}
  Let \Cref{modelassumption:1} and \Cref{step_assumption:inner} hold, and \(f\) be bounded below by $f_{\mathrm{low}}$.
  Let $0 < \eta_1 \leq \eta_2 < 1$ and $0 < \gamma_3 \leq 1 < \gamma_2$ be the parameters of \Cref{alg:R2N}.
  Then,
  \begin{equation}%
    \label{eq:inner-iteration-count-bound}
    j_k(\epsilon_k) \leq \alpha_1 \sigma_{\max, k} \frac{f(x_{k,0}) + \tau_k h(c(x_{k,0})) - f_{\mathrm{low}}}{\epsilon_k^2} + \alpha_2 \left| \log \left(\frac{\sigma_{\max, k}}{\sigma_{k,0}} \right) \right|,
  \end{equation}
  where \(\sigma_{\max, k} = \max(\sigma_{k,0}, \alpha_3 \kappa_{m,k} )\), and
  \begin{equation}%
    \label{eq:def-const-alpha}
    \alpha_1 \coloneq (1 + |\log_{\gamma_1}(\gamma_3)|) > 0, \quad
    \alpha_2 \coloneq | 1 / \log(\gamma_1) | > 0,            \quad
    \alpha_3 \coloneq 2 \gamma_2 / (1 - \eta _2) > 2.
  \end{equation}
\end{theorem}

Because our approach is similar to that of \citep{cartis-gould-toint-2011}, the next lemmas give estimates between our measures, based on $\xi$ and $\theta$, and a generalization of theirs.
For all $x \in \R^n$ and $\tau > 0$, the generalization of their measures \citep[Equations~\((3.1)\) and~\((3.10)\)]{cartis-gould-toint-2011} to any exact norm penalty are
\begin{subequations}%
  \begin{align}%
    \label{eq:xiTR}
    \xi_{\mathrm{TR}}(x, \tau) & \coloneq f(x) + \tau h(c(x)) - \min_{\| s \|_2 \leq 1} (\varphi + \tau \psi)(s;x), \\
    \label{eq:thetaTR}
    \theta_{\mathrm{TR}}(x)    & \coloneq h(c(x)) - \min_{\| s \|_2 \leq 1} \psi(s;x).
  \end{align}
\end{subequations}
Their measures use the \(\ell_2\)-norm for \(\psi\) and \(h\) because their analysis is limited to the $\ell_2$-norm penalty.

\begin{lemma}%
  \label{lemma:xiTR}
  For any $x \in \R^n$, $\tau > 0$ and $\sigma > 0$,
  \begin{equation}%
    \label{eq:xi-xiTR}
    \xi(x; \sigma, \tau) \geq \tfrac{1}{2} \min(1, \, \sigma^{-1} \xi_{\mathrm{TR}}(x, \tau)) \, \xi_{\mathrm{TR}}(x, \tau).
  \end{equation}
  In particular, for any $\epsilon > 0$, if $\sigma^{1/2} \xi{(x ; \sigma, \tau)}^{1/2}  \leq \epsilon$, then $\sigma \geq  \sqrt{2} \epsilon $ implies
  \begin{equation}%
    \label{eq:xiTR-xi}
    \sigma^{1/2} \xi{(x ; \sigma, \tau)}^{1/2} \geq \tfrac{1}{\sqrt{2}}\xi_{\mathrm{TR}}(x, \tau).
  \end{equation}
\end{lemma}

\begin{proof}
  Since $f$ and $c$ are \(\mathcal{C}^1\), and $\tau h$ is convex and globally Lipschitz continuous, \citep[Lemma~\(2.5\)]{cartis-gould-toint-2011} gives
  \(
  \xi(x,\sigma,\tau) - \tfrac{1}{2} \sigma \| s \|^2_2 \geq \tfrac{1}{2}\min (1,\sigma^{-1} \xi_{\mathrm{TR}}(x, \tau) ) \xi_{\mathrm{TR}}(x, \tau)
  \).
  Since $\xi(x,\sigma,\tau) \geq \xi(x,\sigma,\tau) - \tfrac{1}{2} \sigma \| s \|^2_2 $,~\eqref{eq:xi-xiTR} holds.

  For the second part, let \(\epsilon > 0\) and assume that $\sigma^{1/2} \xi{(x ; \sigma, \tau)}^{1/2}  \leq \epsilon$ and \(\sigma \geq \sqrt{2} \epsilon\).
  We now show that \(\sigma^{-1} \xi_{\mathrm{TR}}(x, \tau) \leq 1\), as that will imply~\eqref{eq:xiTR-xi} by way of~\eqref{eq:xi-xiTR}.
  Assume by contradiction that
  \(
  \sigma^{-1}\xi_{\mathrm{TR}}(x,\tau) > 1
  \).
  Our assumption and~\eqref{eq:xi-xiTR} imply
  \(
  \epsilon^2 \sigma^{-1} \geq \xi(x ; \sigma, \tau) \geq \tfrac{1}{2} \xi_{\mathrm{TR}}(x,\tau)
  \).
  We multiply both sides by $\sigma^{-1}$, and obtain
  \(
  \epsilon^2\sigma^{-2} \geq \tfrac{1}{2} \sigma^{-1} \xi_{\mathrm{TR}}(x,\tau) > \tfrac{1}{2}
  \),
  which contradicts our assumption that \(\sigma \geq \sqrt{2} \epsilon\).
\end{proof}

\begin{lemma}%
  \label{lemma:thetaTR}
  For any $x \in \R^n$,
  \(
  \theta_{\mathrm{TR}}(x) \geq \min(\tfrac{1}{\sqrt{2}}, \, \theta{(x)}^{1/2}) \, \theta{(x)}^{1/2}
  \).
\end{lemma}

\begin{proof}
  For $\{ s^* \} = \prox{1\psi}(0) = \argmin{s} \psi(s;x) + \tfrac{1}{2}\| s \|^2_2$, \Cref{prop:prox-decrease} implies
  \begin{equation}%
    \label{eq:thetaTR-1}
    \theta(x) = \psi(0;x) - \psi(s^*;x) \geq \tfrac{1}{2} \| s^* \|^2_2.
  \end{equation}
  Assume first that $\theta(x) \leq \tfrac{1}{2}$.
  Then,~\eqref{eq:thetaTR-1} implies \(\| s^* \|_2 \leq 1\).
  Hence,
  \begin{equation*}%
    \min_{\|s\|_2 \leq 1} h(c(x) + J(x)s) = \min_{\|s\|_2 \leq 1} \psi(s;x) \leq \psi(s^*;x),
  \end{equation*}
  which, together with~\eqref{eq:theta} and~\eqref{eq:thetaTR}, implies that \(\theta_{\mathrm{TR}}(x) \geq \theta(x)\).

  Assume next that $\theta(x) > \tfrac{1}{2}$.
  Note that~\eqref{eq:thetaTR-1} may be written
  \(
  \| \tfrac{1}{\sqrt{2}} \theta{(x)}^{-1/2} s^* \|_2 \leq 1
  \).
  Since $s \to h(c(x) + J(x)s)$ is convex and $\tfrac{1}{\sqrt{2}} \theta{(x)}^{-1/2} < 1$,
  \begin{align*}%
    \min_{\|s\|_2 \leq 1} h(c(x) & + J(x)s) \leq h(c(x) + \tfrac{1}{\sqrt{2}} \theta{(x)}^{-1/2} J(x) s^*)                                                            \\
                                 & \leq \left(1 - \tfrac{1}{\sqrt{2}} \theta{(x)}^{-1/2} \right) h(c(x)) + \tfrac{1}{\sqrt{2}} \theta{(x)}^{-1/2} h(c(x) + J(x) s^*).
  \end{align*}
  We then have from~\eqref{eq:theta} and~\eqref{eq:thetaTR},
  \(
  \theta_{\mathrm{TR}}(x) \geq \tfrac{1}{\sqrt{2}} \theta{(x)}^{-1/2} \theta(x) = \tfrac{1}{\sqrt{2}} \theta{(x)}^{1/2}
  \).
\end{proof}

\begin{theorem}%
  \label{th:thconv}
  Assume that \Cref{step_assumption:kappa-tauk} holds and that \(f\) is bounded below by $f_{\mathrm{low}}$.
  Assume that there is $\bar{\tau} \geq 0$, independent of $\epsilon$, such that $\theta{(x_{k+1})}^{1/2} < \epsilon$ whenever $\tau_k \geq \bar{\tau}$.
  Then \Cref{alg:exactpen} with stopping criterion~\eqref{eq:stop-crit} terminates either with an approximate KKT point of~\eqref{eq:nlp} or with an infeasible critical point of the feasibility measure~\eqref{eq:theta} in at most
  \begin{equation*}%
    \left(\left\lceil \frac{\bar{\tau} - \tau_0}{\beta_1} \right\rceil + \left\lceil \log_{\beta_2} (\epsilon/\epsilon_0) \right\rceil \right)
    \left(
    \frac{\alpha_1 (\kappa_h \bar{\tau} + \kappa_f) (\alpha_4 \bar{\tau} + \alpha_5)}{\epsilon^2} + \alpha_2 \left| \log \left(\frac{\alpha_4 \tau_0 + \alpha_5}{\beta_3 \tau_0} \right) \right|
    \right)
  \end{equation*}
  inner iterations, which is an overall complexity of $\mathcal{O}((\bar{\tau} + |\log(\epsilon/\epsilon_0)|) \, \bar{\tau}^2 \epsilon^{-2})$, where $\alpha_1$, $\alpha_2$, and $\alpha_3$ are defined in~\eqref{eq:def-const-alpha}, $\kappa_{m_1}$ and $\kappa_{m_2}$ are defined in \Cref{step_assumption:kappa-tauk}, and $\kappa_f$, $\kappa_h$ and $\alpha_4$ are defined as
  \begin{subequations}%
    \begin{alignat}{2}%
      \label{eq:kappa-h}
      \kappa_h & \coloneq \max_{k < k(\epsilon)} h(c(x_{k,0})) > 0,                    & \qquad
      \kappa_f & \coloneq \max_{k < k(\epsilon)} f(x_{k,0}) - f_{\mathrm{low}} \geq 0,          \\
      \label{eq:alpha-4}
      \alpha_4 & \coloneq \max(\beta_3, \alpha_3 \kappa_{m_1}) > 0,                    & \qquad
      \alpha_5 & \coloneq \max(\beta_4, \alpha_3 \kappa_{m,2}) > 0.
    \end{alignat}
  \end{subequations}

  Alternatively, assume $\tau_k$ grows unbounded with k.
  Assume further that there is $\kappa_g \geq 0$ such that $\| \nabla f(x_k) \|_2 \leq \kappa_g$ for all \(k\).
  Then, \Cref{alg:exactpen} with stopping criterion~\eqref{eq:stop-crit} terminates either with an approximate KKT point of~\eqref{eq:nlp} or with an infeasible critical point of the feasibility measure~\eqref{eq:theta} in as many iterations as in the first case, but replacing $\bar{\tau}$ with $(\kappa_g + 1) \epsilon^{-2} + 1$, which is an overall complexity of $\mathcal{O}(\epsilon^{-8})$.
\end{theorem}

\begin{proof}
  Consider the first part of the theorem.
  From \Cref{eq:exactpen_epsilon_sequence} of \Cref{alg:exactpen}, there are at most $\lceil \log(\epsilon/\epsilon_0) / \log(\beta_2) \rceil = \lceil \log_{\beta_2} (\epsilon/\epsilon_0) \rceil$ outer iterations where \(\tau_k\) is not increased until ${\theta(x_k)}^{1/2}$ attains or drops below $\epsilon$.
  Similarly, there are at most \(\lceil (\bar{\tau} - \tau_0) / \beta_1 \rceil\) outer iterations where \(\tau_k\) is increased until it reaches or attains \(\bar{\tau}\).
  Hence, $\tau_k$ attains or exceeds $\bar{\tau}$ in
  \begin{equation}%
    \label{eq:k-eps-bound}
    k(\epsilon) \leq \left\lceil \frac{\bar{\tau} - \tau_0}{\beta_1} \right\rceil + \left\lceil \log_{\beta_2} (\epsilon/\epsilon_0) \right\rceil
  \end{equation}
  outer iterations.
  We now find an estimate on $j_k(\epsilon_k)$.
  By \Cref{step_assumption:inner,step_assumption:kappa-tauk} and \Cref{th:inner-complexity},~\eqref{eq:inner-iteration-count-bound} holds.
  First, from~\eqref{eq:kappa-h}
  \begin{equation}%
    \label{eq:fc-bound}
    f(x_{k,0}) + \tau_k h(c(x_{k,0})) - f_{\mathrm{low}} \leq \kappa_h \tau_k + \kappa_f.
  \end{equation}
  The form of \(\sigma_{k, 0}\) in \Cref{eq:exactpen_inner_stationnarity} of \Cref{alg:exactpen}, \Cref{step_assumption:kappa-tauk},~\eqref{eq:alpha-4} and the definition of \(\sigma_{\max, k}\) in \Cref{th:inner-complexity} give
  \begin{equation}%
    \label{eq:nu-min-k-bound}
    \sigma_{\max, k} = \max(\sigma_{k,0} , \alpha_3\kappa_{m,k}) \leq \max(\beta_3 \tau_k, \beta_4, \alpha_3 (\kappa_{m_1}\tau_k + \kappa_{m_2})) \leq \alpha_4 \tau_k + \alpha_5,
  \end{equation}
  so that
  \begin{equation}%
    \label{eq:nu-min-frac-bound}
    \frac{\sigma_{\max, k}}{\sigma_{k, 0}} \leq \frac{\alpha_4 \tau_k + \alpha_5}{\max(\beta_3 \tau_k, \beta_4)} \leq \frac{\alpha_4 \tau_k + \alpha_5}{\beta_3 \tau_k} \leq \frac{\alpha_4 \tau_0 + \alpha_5}{\beta_3 \tau_0},
  \end{equation}
  because $\tau \mapsto \tfrac{\alpha_4 \tau + \alpha_5}{\beta_3 \tau}$ is decreasing on $\R^+$, and $\tau_k \geq \tau_0$.
  Combining~\eqref{eq:fc-bound},~\eqref{eq:nu-min-k-bound} and~\eqref{eq:nu-min-frac-bound} into~\eqref{eq:inner-iteration-count-bound} gives
  \begin{align*}%
    j_k(\epsilon_k) & \leq \frac{\alpha_1 (\kappa_h \tau_k + \kappa_f) (\alpha_4 \tau_k + \alpha_5)}{\epsilon_k^2} + \alpha_2 \left| \log \left(\frac{\alpha_4 \tau_0 + \alpha_5}{\beta_3 \tau_0} \right) \right|
    \\
                    & \leq  \frac{\alpha_1 (\kappa_h \bar{\tau} + \kappa_f) (\alpha_4 \bar{\tau} + \alpha_5)}{\epsilon^2} + \alpha_2 \left| \log \left(\frac{\alpha_4 \tau_0 + \alpha_5}{\beta_3 \tau_0} \right) \right|.
  \end{align*}
  Thus, we obtain the desired bound on \(\sum_{k = 0}^{k(\epsilon)-1} j_k(\epsilon_k)\).

  We now turn to the second part of the theorem.
  First note that whenever $\tau_k$ is increased, \Cref{eq:exactpen_epsilon_sequence} of \Cref{alg:exactpen} implies
  \begin{equation}%
    \label{eq:xi-upper}
    \sigma_{k,j}^{1/2} \xi{(x_k ; \sigma_{k,j},\tau_k)}^{1/2} \leq \epsilon_k < \theta^{1/2}(x_k)
    \quad \text{with} \quad
    j = j_k(\epsilon_k).
  \end{equation}
  Using the notation of \Cref{lemma:xiTR,lemma:thetaTR}, if we assume that $\sigma_{k,j} \geq \sqrt{2}\epsilon_k $, then~\eqref{eq:xi-upper} and \Cref{lemma:xiTR} imply
  \begin{equation}%
    \label{eq:xiTR-xi-theorem}
    \sigma_{k,j}^{1/2} \xi{(x_k ; \sigma_{k,j},\tau_k)}^{1/2} \geq \tfrac{1}{\sqrt{2}} \xi_{\mathrm{TR}}(x_k, \tau_k).
  \end{equation}
  From \Cref{alg:exactpen}, $\sigma_{k,j} \geq \beta_4 > 0$ for all $k$ and $j$.
  Because $\epsilon_k \leq \beta_4 / \sqrt{2}$ for all sufficiently large \(k\), we may assume without loss of generality that~\eqref{eq:xiTR-xi-theorem} holds.
  Furthermore, for any $s_{k,\mathrm{TR}} \in \argmin{\| s \|_2 \leq 1} \psi(s; x_k)$,
  \begin{align}%
    \xi_{\mathrm{TR}}(x_k, \tau_k) & = f(x_k) + \tau_k h(x_k) - \min_{\| s \|_2 \leq 1} (\varphi + \tau_k \psi)(s;x_k)           \\
                                   & \geq f(x_k) + \tau_k h(x_k) - (\varphi + \tau_k \psi)(s_{k,\mathrm{TR}};x_k)                \\
                                   & = \tau_k (h(x_k) - \min_{\| s \|_2 \leq 1} \psi(s;x)) - {\nabla f(x_k)}^T s_{k,\mathrm{TR}} \\
                                   & = \tau_k \theta_{\mathrm{TR}} (x_k) - {\nabla f(x_k)}^T s_{k,\mathrm{TR}}                   \\
                                   &
    \label{eq:xiTR-thetaTR}%
    \geq \tau_k \theta_{\mathrm{TR}} (x_k) - \| \nabla f(x_k) \|_2,
  \end{align}
  where we used the definition of the models~\eqref{eq:models} on the third line, the definition of $\theta_{\mathrm{TR}}$~\eqref{eq:thetaTR} on the fourth and the Cauchy-Schwarz inequality on the last line combined with the fact that $\|s_{k,\mathrm{TR}}\|_2 \leq 1$.
  Substituting~\eqref{eq:xiTR-thetaTR} into~\eqref{eq:xiTR-xi-theorem} then gives
  \begin{equation}%
    \label{eq:xi-thetaTR}
    \sigma_{k,j}^{1/2} \xi{(x_k ; \sigma_{k,j},\tau_k)}^{1/2} \geq \tfrac{1}{\sqrt{2}} \left(\tau_k \theta_{\mathrm{TR}} (x_k) - \| \nabla f(x_k) \| \right).
  \end{equation}

  At this point, we consider two cases.
  Assume first that ${\theta(x_k)}^{1/2} \leq \tfrac{1}{\sqrt{2}}$.
  \Cref{lemma:thetaTR} and~\eqref{eq:xi-thetaTR} combine with~\eqref{eq:xi-upper} to give
  \begin{equation}%
    \label{eq:theta-theta}
    \tfrac{1}{\sqrt{2}} \left(\tau_k \theta (x_k) - \| \nabla f(x_k) \|_2 \right) \leq \sigma_{k,j}^{1/2} {\xi(x_k ; \sigma_{k,j},\tau_k)}^{1/2} \leq {\theta(x_k)}^{1/2} \leq \tfrac{1}{\sqrt{2}} (\theta(x_k) + 1),
  \end{equation}
  whenever $\tau_k$ is increased, where we used the fact that $\sqrt{t} \leq \tfrac{t+1}{2} \leq \tfrac{t+1}{\sqrt{2}}$ for any $t > 0$ in the last inequality.
  Thus,~\eqref{eq:theta-theta} implies
  \begin{equation}%
    \label{eq:theta-theta1}
    \| \nabla f(x_k) \|_2 \geq (\tau_k - 1) \theta(x_k) - 1.
  \end{equation}
  Because $\tau_k > 1$ for all large enough $k$, and $\| \nabla f(x_k) \|_2 \leq \kappa_g$,~\eqref{eq:theta-theta1} yields
  \(
  \theta(x_k) \leq (\kappa_g + 1) / (\tau_k - 1)
  \)
  whenever \(\tau_k\) is increased.
  Therefore, $\theta^{1/2}(x_k) \leq \epsilon_k$ whenever
  \begin{equation}%
    \label{eq:tauk-lower}
    \tau_k \geq \frac{\kappa_g + 1}{\epsilon_k^2} + 1.
  \end{equation}

  Assume instead that $\theta^{1/2}(x_k) > \tfrac{1}{\sqrt{2}}$.
  \Cref{lemma:thetaTR} and~\eqref{eq:xi-thetaTR} imply
  \begin{equation}%
    \label{eq:xi-theta1}
    \sigma_{k,j}^{1/2} \xi{(x_k ; \sigma_{k,j},\tau_k)}^{1/2} \geq \tfrac{1}{\sqrt{2}} \left(\tfrac{1}{\sqrt{2}} \tau_k {\theta(x_k)}^{1/2} - \| \nabla f(x_k) \|_2 \right).
  \end{equation}
  Using the same approach that led from~\eqref{eq:theta-theta} to~\eqref{eq:tauk-lower}, with the difference that we do not require the final inequality in~\eqref{eq:theta-theta}, we find in this case
  \begin{equation}%
    \label{eq:tauk-lower_bis}
    \tau_k \geq \frac{\sqrt{2} \kappa_g}{\epsilon_k} + 2.
  \end{equation}
  Combining~\eqref{eq:tauk-lower} with~\eqref{eq:tauk-lower_bis} gives ${\theta(x_k)}^{1/2} \leq \epsilon_k $ whenever
  \begin{equation}%
    \label{eq:tauk-lower-max}
    \tau_k \geq \max \left( \frac{\kappa_g + 1}{\epsilon_k^2} + 1, \, \frac{\sqrt{2} \kappa_g}{\epsilon_k} + 2 \right).
  \end{equation}
  For small enough values of $\epsilon_k$, i.e., for all large enough \(k\), $\tfrac{\kappa_g + 1}{\epsilon_k^2} + 1 > \tfrac{\sqrt{2} \kappa_g}{\epsilon_k} + 2$, and we may assume without loss of generality that ${\theta(x_k)}^{1/2} \leq \epsilon_k $ whenever~\eqref{eq:tauk-lower} holds.

  Overall,~\eqref{eq:tauk-lower} implies $\theta^{1/2}(x_k) \leq \epsilon$ whenever
  \begin{equation}%
    \label{eq:tauk-lower1}
    \tau_k \geq \frac{\kappa_g + 1}{\epsilon^2} + 1.
  \end{equation}
  The value on the right-hand side of~\eqref{eq:tauk-lower1} is the counterpart of $\bar{\tau}$ from the first part of the proof, with the crucial difference that it depends on $\epsilon$.
  Still, we may use the value of the right-hand side of~\eqref{eq:tauk-lower1} in place of $\bar{\tau}$ in the first part of the proof.
\end{proof}

The existence of $\kappa_h$ and $\kappa_f$ in~\eqref{eq:kappa-h} is guaranteed by the fact that $f$, $h$ and $c$ are continuous functions.

Combining \Cref{th:exact-penanlty} with \Cref{th:thconv}, we see that the assumptions of the first case of \Cref{th:thconv} hold under the MFCQ\@.
Therefore, the bound in the first case of \Cref{th:thconv} is the one we expect to observe the most in practical settings, while the bound in the second case should only happen if the problem is degenerate, i.e., does not admit Lagrange multipliers.

Comparing the bounds in \Cref{th:thconv} with those of \citet[Theorem~\(3.2\)]{cartis-gould-toint-2011}, we see two differences: an additional $\log(\epsilon/\epsilon_0)$ term in our first case, and we found an $\mathcal{O}(\epsilon^{-8})$ instead of their $\mathcal{O}(\epsilon^{-5})$ complexity.

The additional $\log(\epsilon/\epsilon_0)$ comes from the fact that we used an increasing accuracy strategy in \Cref{alg:exactpen}, as we do not need to compute high-accuracy solutions when far from feasibility.
However, if one chooses to enforce constant accuracy throughout, as is done in \citep{cartis-gould-toint-2011}, then $\epsilon_0 = \epsilon$, which implies $\log(\epsilon/\epsilon_0) = 0$.

Regarding the difference in the second case, notice first that our measure $\theta^{1/2}$ does not scale linearly with \citeauthor{cartis-gould-toint-2011}’s $\theta_{\mathrm{TR}}$ but $\theta$ does.
We can easily find examples for which \(\theta(x) = \theta_{\mathrm{TR}}(x)\).
Consider $c : \R \to \R$, $c(x) = x$.
For \(x \in \R\), from~\eqref{eq:modpsi}, \(\psi(s; x) = |x + s|\) and
\begin{equation*}%
  \{s^*\} = \argmin{|s| \leq 1} \psi(s; x) = \prox{1\psi}(0) =
  \begin{cases}
    \{ -1 \}           & \text{if } x > 1,      \\
    \{ -x \}           & \text{if } |x| \leq 1, \\
    \{ \phantom{-}1 \} & \text{if } x < -1,
  \end{cases}
\end{equation*}
which implies from~\eqref{eq:theta} and~\eqref{eq:thetaTR} that
\begin{align*}%
  \theta(x)               & = |x| - |x + s^*| = \min(|x|, \, 1), \\
  \theta_{\mathrm{TR}}(x) & = |x| - |x + s^*| = \min(|x|, \, 1).
\end{align*}
Now, we can see that the issue really comes from \Cref{lemma:thetaTR}; if instead of the inequality in \Cref{lemma:thetaTR} we had an inequality of the type $\theta_{\mathrm{TR}} > \kappa\theta^{1/2}$ for some constant $\kappa > 0$, then, using the same approach that led to~\eqref{eq:tauk-lower_bis}, we would have found an upper bound on $\tau_k$ in $\mathcal{O}(\epsilon^{-1})$ instead of $\mathcal{O}(\epsilon^{-2})$ which would have then led to the $\mathcal{O}(\epsilon^{-5})$ complexity.

Furthermore, if we consider now \(c : \R \to \R\), \(c(x) = x^2/8 \).
For \(x \in \R\), from~\eqref{eq:modpsi}, \( \psi(s; x) = \tfrac{1}{4} | \tfrac{1}{2}x^2 + xs | \) and
\begin{align*}%
  \{s^*\} = \prox{1\psi}(0) = & \argmin{s} \tfrac{1}{2} s^2 + \psi(s; x) = \{ -\tfrac{1}{2} x \} \\
  \{s^*_{\mathrm{TR}}\} =     & \argmin{|s| \leq 1} \psi(s; x) =
  \begin{cases}
    \{ -1 \}            & \text{if } x > 2,      \\
    \{ -\tfrac{x}{2} \} & \text{if } |x| \leq 2, \\
    \{ \phantom{-}1 \}  & \text{if } x < -2.
  \end{cases}
\end{align*}
Then, from~\eqref{eq:theta} and~\eqref{eq:thetaTR},
\begin{align*}%
  \theta(x)               & = \tfrac{1}{8}|x^2| - \tfrac{1}{4}|\tfrac{1}{2} x^2 + x s^*| = \tfrac{1}{8}x^2.                                         \\
  \theta_{\mathrm{TR}}(x) & = \tfrac{1}{8}|x^2| - \tfrac{1}{4}|\tfrac{1}{2} x^2 + x s^*_{\mathrm{TR}}| = \tfrac{1}{4}\min(|x|, \, \tfrac{1}{2}x^2).
\end{align*}
Therefore in this case, for all \(x \in \R\), \(\theta_{\mathrm{TR}}(x) = \min(\tfrac{1}{\sqrt{2}}, \, \theta{(x)}^{1/2})\theta{(x)}^{1/2}\) so that the bound in \Cref{lemma:thetaTR} is tight.

For $\{ s^* \} = \prox{1\psi}(0) = \arg\min_{s} \tfrac{1}{2}\|s\|^2_2 + \psi(s;x)$,
\begin{equation}
  \label{eq:psi-subgradient}
  -s^* \in \partial \psi(s^*;x).
\end{equation}
For any $g_{s^*,x} \in \partial \psi(s^*;x)$, convexity of $\psi(\cdot;x)$ and \citep[Corollary D.2.1.3]{hiriart-urruty-lemaréchal-2001} imply
\begin{equation*}
  \psi(s^*;x) = \psi(0;x) + g_{s^*,x}^T s^* + o(\|s^*\|_2).
\end{equation*}
The choice $g_{s^*, x} = -s^*$,~\eqref{eq:theta} and~\eqref{eq:psi-subgradient} combine to give
\(
\theta(x) =  \|s^*\|^2_2 + o(\|s^*\|_2)
\).
Hence, $\theta$ varies like the squared norm of a subgradient, which justifies taking the square root of $\theta$ to define a stationarity measure.
In addition, \Cref{lemma:thetaTR} implies $\theta_{\mathrm{TR}} \geq \theta$ as $\theta \to 0$.
The above leads us to believe that the $\mathcal{O}(\epsilon^{-8})$ is a more realistic worst-case complexity bound for exact penalty algorithms when the penalty parameter is unbounded, though neither bound has been proven to be tight.

Finally, the assumptions of \Cref{th:thconv} are weaker than those of \citet{cartis-gould-toint-2011} as we only require boundedness of $\{ \nabla f(x_k) \}$, whereas they require boundedness of $\{ x_k \}$.

\subsection*{A quasi-Newton variant}%
\label{sec:qN}

Instead of using \citep[Algorithm~\(6.1\)]{aravkin-baraldi-orban-2022} as solver for~\eqref{eq:penalty-nlp} in \Cref{alg:exactpen}, we explore variant R2N recently proposed in \citep{diouane-habiboullah-orban-2024} that uses or approximates second-order information on $f$, and show that it preserves convergence and complexity properties.
As it turns out, using second-order models in the subproblem does not improve the complexity bounds of the outer loop, but we still expect improvements in practice as more information is given to enrich the model.
At inner iteration \(j\) of outer iteration \(k\), instead of~\eqref{eq:prox-problem}, we compute \(s_{k, j}\) as an approximate solution of
\begin{equation}%
  \label{eq:prox-problem-quadratic}
  \minimize{s} \ m_Q(s; x_{k, j}, B_{k, j}) + \tfrac{1}{2}\sigma_{k, j} \|s\|^2_2,
\end{equation}
where \(m_Q(s; x, B) \coloneqq \varphi_Q(s ; x, B) + \tau_k \psi(s ; x)\), and
\(
\varphi_Q(s; x, B) \coloneqq \varphi(s; x) + \tfrac{1}{2} s^T B s
\),
with $\varphi$ defined in~\eqref{eq:modphi1}, $B = B^T \in \R^{n \times n}$.
In particular, $B$ may be a quasi-Newton approximation of $\nabla^2 f(x)$, or $\nabla^2 f(x)$ if it exists.
\citet{diouane-habiboullah-orban-2024} do not require \(B\) to be positive (semi-)definite.
However, it is not difficult to see that~\eqref{eq:prox-problem-quadratic} is unbounded below along any negative curvature direction of \(B + \sigma I\).
Hence, we always adjust \(\sigma\) so that \(\varphi_Q(s; x, B) + \tfrac{1}{2} \sigma \|s\|^2_2\) is convex.
Clearly, \Cref{modelassumption:1} continues to hold for \(\varphi_Q\).

Although the analysis of \citep{diouane-habiboullah-orban-2024} does not depend on it, we make the following simplifying assumption.

\begin{stepassumption}%
  \label{step_assumption:inner_qn}
  There is a constant $\kappa_B$ such that for each iteration $k$ and $j$ of \Cref{alg:exactpen} that uses R2N,
  \(
  \| B_{k, j} \|_2 \leq \kappa_B
  \).
\end{stepassumption}

R2N uses $\xi(x_{k, j}; \nu_{k, j}^{-1}, \tau_k)$, defined in~\eqref{eq:xi}, as a stationarity measure, where
\(
\nu_{k, j}^{-1} \coloneqq \theta_1 / (\sigma_{k, j} + \| B_{k, j} \|_2)
\),
and where $\theta_1 \in (0, \, 1)$ is a parameter.

An analogue of \Cref{th:thconv} holds for \Cref{alg:exactpen} with R2N as subsolver.
If we denote $\tilde\sigma_{k, j} \coloneqq {\nu_{k, j}}^{-1}$, \Cref{step_assumption:inner} yields
\begin{equation*}
  \sigma_{k, j}\theta_1^{-1} \leq \tilde\sigma_{k, j} \leq \sigma_{k, j}\theta_1^{-1} + \kappa_B\theta_1^{-1}.
\end{equation*}
With these inequalities, we find $\tilde\sigma_{\max, k} = \theta_1^{-1}\sigma_{\max, k} + \theta_1^{-1}\kappa_B$ which acts as the $\sigma_{\max, k}$ defined previously.
With that, we can establish \Cref{th:thconv} with $\xi(x_{k, j}; \nu_{k, j}^{-1}, \tau_k)$ as inner stationarity measure.

The updated complexity bounds rest upon the following result.

\begin{theorem}[\protect{\citealp[Theorems~\(6.4\) and~\(6.5\)]{diouane-habiboullah-orban-2024}}]%
  \label{th:inner-complexity-R2N}
  Let \Cref{modelassumption:1,step_assumption:inner,step_assumption:inner_qn} hold and \(f\) be bounded below by $f_{\mathrm{low}}$.
  Let $0 < \eta_1 \leq \eta_2 < 1$ and $0 < \gamma_3 \leq 1 < \gamma_2$ be the parameters of \Cref{alg:R2N}.
  Then,
  \begin{equation*}%
    j_k(\epsilon_k) \leq (2\kappa_B( 1 + \tilde\sigma_{\max, k} ) + \tilde\sigma_{\max, k}) \left(1 + |\log_{\gamma_1}(\gamma_3)| \right) \frac{\Delta (f+h)}{\eta_1 \alpha \epsilon_k^2} + \log_{\gamma_1}\left(\frac{\tilde\sigma_{\max, k}}{\sigma_{k,0}} \right),
  \end{equation*}
  where \(\Delta (f+h) \coloneqq f(x_{k,0}) + \tau_k h(c(x_{k,0})) - f_{\mathrm{low}}\).
\end{theorem}

\begin{theorem}%
  \label{th:thconv-R2N}
  Assume that \Cref{step_assumption:kappa-tauk} and the assumptions of \Cref{th:inner-complexity-R2N} hold.
  Assume that there is $\bar{\tau} \geq 0$, independent of $\epsilon$, such that $\theta{(x_{k+1})}^{1/2} < \epsilon$ whenever $\tau_k \geq \bar{\tau}$.
  Then \Cref{alg:exactpen} with stopping criterion~\eqref{eq:stop-crit} using \Cref{alg:R2N} to solve~\eqref{eq:penalty-nlp} terminates either with an approximate KKT point of~\eqref{eq:nlp} or with an infeasible critical point of the feasibility measure~\eqref{eq:theta} in at most
  \begin{equation*}%
    \left(\left\lceil \frac{\bar{\tau} - \tau_0}{\beta_1} \right\rceil + \left\lceil \log_{\beta_2} (\epsilon/\epsilon_0) \right\rceil \right)
    \left(
    \frac{ \alpha_6 }{ \alpha \epsilon^2 } + \alpha_2 \left| \log \left(\frac{\alpha_4 \tau_0 + \alpha_5}{\beta_3 \tau_0} \right) \right|
    \right)
  \end{equation*}
  inner iterations, which is an overall complexity of $\mathcal{O}((\bar{\tau} + |\log(\epsilon/\epsilon_0)|) \, \bar{\tau}^2 \epsilon^{-2})$, where $\alpha_1$, $\alpha_2$, $\alpha_3$ are as in \Cref{th:inner-complexity}, $\alpha_4$, $\alpha_5$, $\kappa_h$, $\kappa_f$ are as in \Cref{th:thconv}, $\kappa_{m_1}$ and $\kappa_{m_2}$ are as in \Cref{step_assumption:kappa-tauk}, $\kappa_B$ is as in \Cref{step_assumption:inner_qn}, and
  \[
    \alpha_6 \coloneqq \alpha_1 (\kappa_h \bar{\tau} + \kappa_f) ((\alpha_4 \bar{\tau} + \alpha_5)(1 + 2\kappa_B ) + 2 \kappa_B).
  \]

  Alternatively, assume $\tau_k$ grows unbounded with k.
  Assume further that there is $\kappa_g \geq 0$ such that $\| \nabla f(x_k) \|_2 \leq \kappa_g$ for all \(k\).
  Then, \Cref{alg:exactpen} with stopping criterion~\eqref{eq:stop-crit} using R2N to solve~\eqref{eq:penalty-nlp}  terminates either with an approximate KKT point of~\eqref{eq:nlp} or with an infeasible critical point of the feasibility measure~\eqref{eq:theta} in as many iterations as in the first case, but replacing $\bar{\tau}$ with $(\kappa_g + 1) \epsilon^{-2} + 1$, which is an overall complexity of $\mathcal{O}(\epsilon^{-8})$.
\end{theorem}

The proof of \Cref{th:thconv-R2N} is almost identical to that of \Cref{th:thconv}.
The only difference is that the bound~\eqref{eq:nu-min-frac-bound} becomes
\(
\tilde\sigma_{\max, k} / \tilde\sigma_{k, 0} \leq (\sigma_{\max, k} + \kappa_B) / \sigma_{k, 0}
\).

\section{Proximal operators}%
\label{sec:prox-operators}

We now provide an explicit characterization of the solution \(s_{k, j, \mathrm{cp}}\) in~\eqref{eq:prox-problem}.
We begin with general results and then specialize them to our setting.
Moreover, when \Cref{alg:R2N} is used, we derive an efficient procedure for computing \(s_{k, j}\) in~\eqref{eq:prox-problem-quadratic} that avoids the need for an iterative solver for nonsmooth problems.
Although the computation may involve solving a scalar equation, we believe it remains more efficient than general-purpose proximal evaluations.
\begin{theorem}%
  \label{th::thProx-dual}
  Let $Q \in \R^{n \times n}$ be symmetric positive definite, $A \in \R^{m \times n}$, $b \in \R^m$, $d \in \R^n$ and $\tau > 0$.
  The unique solution of
  \begin{equation}%
    \label{eq::thProx-dual0}
    \minimize{u \in \R^n} \ \tfrac{1}{2} u^\top Qu - d^\top u + \tau \| Au + b\|,
  \end{equation}
  is
  \begin{equation}%
    \label{eq::thProx-dual-u-y}
    u^* = Q^{-1}(d + A^T y^*)
  \end{equation}
  where $y^*$ is a solution of
  \begin{equation}%
    \label{eq::thProx-dual1}
    \underset{{y \in \R^m}}{\textup{maximize}} \ -\tfrac{1}{2} {(d + A^T y)}^T Q^{-1} (d + A^T y) - b^T y \quad \st \ \| y \|_* \leq \tau.
  \end{equation}
  If $A$ has full row rank, $y^*$ is unique.
\end{theorem}

\begin{proof}%
  As~\eqref{eq::thProx-dual0} is strictly convex, it has a unique solution $u^*$.
  It may be written
  \begin{equation}%
    \label{eq::thProx-dual2}
    \minimize{u \in \R^n \ z \in \R^m} \ \tfrac{1}{2} u^T Qu - d^T u + \tau \|z\| \quad \st \ z = Au +b.
  \end{equation}
  The Lagrangian of~\eqref{eq::thProx-dual2} is
  \(
  \mathcal{L} (u , z ; y) 
  = \left[ \tfrac{1}{2} u^T Qu - d^T u - {(A^T y)}^T u \right] + \left[ \tau \|z\| + y^T z \right] - b^T y
  \),
  and is separable with respect to $u$ and $z$.
  The minimizer of the terms in $u$ is
  \(
  u^* = Q^{-1} (d + A^T y)
  \),
  and corresponds to an optimal value of
  \(
  -\tfrac{1}{2} {(d + A^T y)}^T Q^{-1} (d + A^T y)
  \).
  Regarding the terms in z,
  \begin{equation*}%
    \min_{z} \left[ \tau \|z\| + y^T z \right] = - \max_{z} \left[ {(-y)}^T z - \tau \|z\| \right] = -{(\tau h)}^*(-y),
  \end{equation*}
  where \({(\tau h)}^*\) is the Fenchel conjugate of \(\tau h\). Using \citep[Example 3.26]{boyd-vandenberghe-2004}, ${(\tau h)}^*(w)=\tau \chi(w / \tau \mid \B_*)$.
  We conclude that
  \begin{equation*}%
    \min_{z} \left[ \tau \|z\| + y^T z \right] =
    \begin{cases}
      0       & \| y \|_* \leq \tau, \\
      -\infty & \| y \|_* > \tau.
    \end{cases}
  \end{equation*}
  The above gives us~\eqref{eq::thProx-dual1} as the dual of~\eqref{eq::thProx-dual2}.
  Since $(u,z)=(0,b)$ is feasible for~\eqref{eq::thProx-dual2}, the relaxed Slater condition, hence strong duality, holds for~\eqref{eq::thProx-dual2} and~\eqref{eq::thProx-dual1}.
  If $A$ has full row rank,~\eqref{eq::thProx-dual1} is strictly concave, and therefore has a unique solution.
\end{proof}
In view of \Cref{th::thProx-dual}, computing a step under a general norm penalty reduces to solving a trust-region subproblem in the dual norm.
When \(h = \|\cdot\|_2\), this subproblem can be solved efficiently, as it has been extensively studied over the past decades \citep[\S7]{conn-gould-toint-2000}.
For this reason, we restrict our attention to this special case for the remainder of the paper.

\begin{theorem}%
  \label{th::thProx-dual-p2}
  Under the assumptions of \Cref{th::thProx-dual}, in the case where $\|\cdot\| = \|\cdot\|_2$, a solution of~\eqref{eq::thProx-dual1} is given by
  \begin{equation*}%
    y^* =
    \begin{cases}%
      - {(A Q^{-1} A^T)}^\dagger (AQ^{-1}d +b) \quad           & \textup{if} \ \| {(A Q^{-1} A^T)}^\dagger (AQ^{-1}d +b)\|_2 \leq \tau \\
                                                               & \textup{and} \ AQ^{-1}d + b \in \textup{Range}(A Q^{-1} A^T)          \\
      - {(A Q^{-1} A^T + \alpha^* I)}^{-1} (AQ^{-1}d +b) \quad & \textup{otherwise},
    \end{cases}
  \end{equation*}
  where $\alpha^*$ is the unique positive root of the strictly decreasing function
  \begin{equation*}%
    g(\alpha) = \|{(AQ^{-1}A^T + \alpha I)}^{-1} (AQ^{-1}d + b)\|_2^2 -\tau^2.
  \end{equation*}
\end{theorem}

\begin{proof}%
  In the case where $p = 2$, we can rewrite~\eqref{eq::thProx-dual1} as
  \begin{equation}%
    \label{eq::thProx-dual-p2-1}
    \minimize{y \in \R^m} \ \tfrac{1}{2} {(d + A^T y)}^T Q^{-1} (d + A^T y) + b^T y \quad \st \ \| y \|_2^2 \leq \tau^2.
  \end{equation}
  Whether \(A\) has full row rank or not,~\eqref{eq::thProx-dual-p2-1} is convex and satisfies Slater's condition.
  Thus, by the KKT conditions, $y$ solves~\eqref{eq::thProx-dual-p2-1} if and only if there is $\alpha^*$ such that
  \begin{equation}%
    \label{eq::thProx-dual-p2-2}
    (AQ^{-1}%
    A^T + \alpha^*I)
    y +AQ^{-1}d + b = 0,
    \quad \text{and} \quad
    0 \leq \alpha^* \perp (\tau - \|y\|_2) \geq 0.
  \end{equation}
  There are now two cases.
  In the first case, $\alpha^* = 0$.
  Then by~\eqref{eq::thProx-dual-p2-2},
  \(
  {(A Q^{-1} A^T)}y = - (AQ^{-1}d + b)
  \).
  Thus, if $AQ^{-1}d + b \in \textup{Range}(A Q^{-1} A^T)$, by primal feasibility, the pseudo-inverse is the solution whenever $\|  {(AQ^{-1}A^\top)}^\dagger (AQ^{-1}d + b) \|_2 \leq \tau$.
  In the second case, $\alpha^* > 0$ and~\eqref{eq::thProx-dual-p2-2} yields
  \(
  \| y \|_2^2 = \tau^2
  \)
  and
  \(
  y = - {(AQ^{-1}A^T + \alpha^*I)}^{-1} (AQ^{-1}d + b)
  \),
  because $AQ^{-1}A^\top + \alpha^* I$ is nonsingular, and where $g(\alpha^*) = 0$.
\end{proof}

When \(A\) has full row rank, the solution given by \Cref{th::thProx-dual-p2} is unique.
When \(A\) is rank deficient in \Cref{th::thProx-dual-p2} and \(\alpha^* = 0\), the pseudo-inverse yields the minimum Euclidean norm solution of the linear system in~\eqref{eq::thProx-dual-p2-2}.
If the latter does not have norm at most \(\tau\), no other solution does, and we must be in the case \(\alpha^* > 0\).
The next result specializes \Cref{th::thProx-dual-p2} to the evaluations of the proximal operator for \(\tau\psi(\cdot; x)\).

\begin{corollary}%
  \label{cor:prox}
  Let $A \in \R^{m \times n}$, $b \in \R^m$, and $\tau > 0$.
  Let \(\eta : u \to \tau\|A u + b\|_2\).
  For \(\nu > 0\),
  \begin{equation}%
    \label{eq:prox-value-rank-deficient}
    \prox{\nu \eta}(w) \ni
    \begin{cases}%
      w - A^T{(A A^T)}^\dagger (Aw + b) \quad          & \textup{if} \ \| {(A A^T)}^\dagger (Aw + b)\|_2 \leq \nu\tau \\
                                                       & \textup{and} \ Aw + b \in \textup{Range}(A A^T)              \\
      w - A^T{(A A^T + \alpha^* I)}^{-1} (Aw +b) \quad & \textup{otherwise},
    \end{cases}
  \end{equation}
  where $\alpha^*$ is the unique positive root of the strictly decreasing function
  \begin{equation*}%
    g(\alpha) = \|{(AA^T + \alpha I)}^{-1} (Aw + b)\|_2^2 -\nu^2\tau^2.
  \end{equation*}
\end{corollary}

\begin{proof}%
  By~\eqref{def:moreau-prox},
  \(
  \prox{\nu \eta}(w) 
  = \argmin{u \in \R^n} \tfrac{1}{2} u^T Iu - w^T u+ \nu \tau \|Au + b\|_2
  \).
  Replacing $Q$ with $I_n$, $d$ with $w$ and $\tau$ with $\nu\tau$ in \Cref{th::thProx-dual,th::thProx-dual-p2} concludes.
\end{proof}

In the context of \Cref{alg:R2N}, we will apply \Cref{cor:prox} with \(\tau = \tau_k\), \(\eta := \tau_k \psi(\cdot; x_k)\), \(A = J(x_k)\), \(b = c(x_k)\), and either \(\nu = \sigma_k^{-1}\) or \(\nu = \tfrac{\theta_1}{\|B_k\|_2+\sigma_k}\) depending on whether we use R2 or R2N respectively.

The same methodology allows us to derive an explicit characterization of solutions of~\eqref{eq:prox-problem-quadratic}.

\begin{corollary}%
  \label{cor:prox-quadratic}
  Let $A \in \R^{m \times n}$ have full row rank, let $b \in \R^m$, $x \in \R^n$ and $\tau > 0$.
  Let $B$ be an $n \times n$ symmetric matrix such that $Q := \nu B + I$ is positive definite, and $\eta : u \xrightarrow{} \tau \|Au + b\|_2$.
  For $\nu > 0$, define
  \begin{align}%
    \label{def:prox-quadratic}
    \prox{\nu \eta}(w, B) & \coloneqq \argmin{u \in \R^n} \ \tfrac{1}{2} \nu^{-1} \|u - w\|^2_2 + \eta(u) + \tfrac{1}{2}  u^T B u.
  \end{align}
  Let \(v \coloneqq A Q^{-1} w + b\).
  Then, the only element in \(\prox{\nu \eta}(w, B)\) is
  \begin{equation*}
    \begin{cases}%
      Q^{-1} \left(w - A^T{(A Q^{-1} A^T)}^\dagger v \right)            & \textup{ if } \| {(A Q^{-1} A^T)}^\dagger v\|_2 \leq \nu\tau \\
                                                                        & \textup{ and } v \in \textup{Range}(A Q^{-1} A^T)            \\
      Q^{-1} \left( w - A^T{(A Q^{-1} A^T + \alpha^* I)}^{-1} v \right) & \textup{ otherwise},
    \end{cases}
  \end{equation*}
  where $\alpha^*$ is the unique positive root of the strictly decreasing function
  \begin{equation*}%
    g(\alpha) = \|{(A Q^{-1} A^T + \alpha I)}^{-1} (AQ^{-1} w + b)\|_2^2 -\nu^2\tau^2.
  \end{equation*}
\end{corollary}

\begin{proof}%
  By~\eqref{def:prox-quadratic},
  \(
  \prox{\nu \eta}(w, B) 
  = \argmin{u \in \R^n} \tfrac{1}{2} u^T Q u - w^T u + \nu \tau \|Au + b\|_2
  \).
  Replacing $Q$ with $\nu B + I$, $d$ with $w$ and $\tau$ with $\nu\tau$ in \Cref{th::thProx-dual,th::thProx-dual-p2} concludes.
\end{proof}

An efficient procedure to solve the $\ell_2$-norm trust-region subproblem was proposed by \citet{more-sorensen-1983}.
We follow their approach.
In the light of \Cref{cor:prox} and using the same notation, it should be clear that any efficient procedure to solve~\eqref{eq::thProx-dual1} should try to find the root of $g(\alpha)$ in as few iterations as possible because each evaluation of $g$ requires solving a linear system.
The idea in \citep{more-sorensen-1983} is to solve the equivalent \emph{secular} equation
\begin{equation}
  \label{eq::secular-equation}
  g(\alpha) = 0 \quad \Longleftrightarrow \quad \frac{1}{\|{(A A^T + \alpha I)}^{-1} (Ax + b)\|_2} = \frac{1}{\nu \tau},
\end{equation}
which has preferable numerical properties.
They use Newton’s method to solve~\eqref{eq::secular-equation}.
We summarize the procedure as \Cref{alg:proximal} with the help of the following result.

\begin{lemma}[\protect{\citealp[Lemma 7.3.1]{conn-gould-toint-2000}}]%
  \label{lemma:secular-derivative}
  Define \(\displaystyle \phi(\alpha) \coloneqq \frac{1}{\|s(\alpha)\|_2} - \frac{1}{\nu \tau}\), where
  \[
    s(\alpha) \coloneqq
    \begin{cases}
      - {(AA^T)}^{\dagger} (Ax + b)       & \quad \text{if } \alpha = 0, \\
      - {(AA^T + \alpha I)}^{-1} (Ax + b) & \quad \text{if } \alpha > 0.
    \end{cases}
  \]
  Then, \(\phi\) is strictly increasing and concave, and for all \(\alpha > 0\),
  \begin{equation*}
    \phi'(\alpha) = - \frac{ {s(\alpha)}^T \nabla_\alpha s(\alpha)}{\| s(\alpha) \|^3_2},
    \quad
    \nabla s(\alpha) =  - {(AA^T + \alpha I)}^{-1} s(\alpha).
  \end{equation*}
  If \(A\) has full row rank then the above holds for $\alpha = 0$ as well.
\end{lemma}

\begin{algorithm}[ht]%
  \caption{%
    \label{alg:proximal}
    Proximal operator evaluation.
  }
  \begin{algorithmic}[1]%
    \State%
    \label{eq:s(0)}
    Compute $s(0)$ as given by \Cref{lemma:secular-derivative}.

    \State%
    If $\|s(0)\| \leq \nu\tau$ and \(A A^T s(0) - (Ax + b) = 0\), return $x + A^T s(0)$

    \State%
    Choose $\theta \in (0,1)$.
    Set \(k = 0\).

    \State%
    \label{eq:alg_proximal_init}
    If $A$ has full row rank, set $\alpha_0 \coloneqq 0$.
    Otherwise, set $\alpha_0 \coloneqq \sqrt{\epsilon_M}$ and compute $s(\alpha_0)$.

    \While{$\| s(\alpha_k) \| \neq \nu\tau$}%

    \State%
    Compute $\alpha_{k+1}$ as $\alpha^+$ given by~\eqref{eq:newton-update}.
    \State%
    \label{eq:newton-safeguard}
    If $\alpha^+ \leq 0$, reset $\alpha_{k + 1} \coloneqq \theta \alpha_k$.
    \State%
    Compute $s(\alpha_{k+1})$.

    \State%
    Let $k \coloneqq k+1$.
    \EndWhile%

    \State{\Return{$ x + A^T s(\alpha_k) $}.}%
  \end{algorithmic}
\end{algorithm}

Instead of computing the Cholesky factorization of $AA^T + \alpha I$, we compute the $QR$ factorization of $A_{\alpha}^T \coloneqq
  \begin{pmatrix}%
    A & \sqrt{\alpha}I_m
  \end{pmatrix}
  ^T$.
The \(R\) factor coincides with the Cholesky factor.
The update of Newton's method can be found using \Cref{lemma:secular-derivative}:
\begin{equation}
  \label{eq:newton-update}
  \alpha^+ = \alpha_k - \frac{\phi(\alpha_k)}{\phi'(\alpha_k)} = \alpha_k - \frac{(\nu \tau - \|s(\alpha_k)\|) \, \|s(\alpha_k)\|^2}{\nu \tau s{(\alpha_k)}^T \nabla s(\alpha_k)}.
\end{equation}
\Cref{alg:proximal} obtains \(q_k = s(\alpha_k)\) by solving \(R^T R q_k = -(Ax + b)\).
Instead of performing another solve with \(R\) and \(R^T\) to obtain \(\nabla s(\alpha_k)\), we observe that \(s{(\alpha_k)}^T \nabla s(\alpha_k) = -q_k^T R^{-1} R^{-T} q_k = -\|R^{-T} q_k\|^2 = -\|p_k\|^2\), so only one extra solve with \(R^T\) is necessary.

When evaluating $s(0)$ on \Cref{eq:s(0)}, we check whether $A$ is rank deficient or not by looking for zeros on the diagonal of $R$.
If so, we compute the least-norm solution of the underdetermined linear system.
\citep[Lemma 7.3.2]{conn-gould-toint-2000} indicates that convergence is assured once \(\alpha \in \mathcal{L} \coloneqq (\max(0, \, -\lambda_1), \, \alpha^*]\), where $\alpha^*$ solves~\eqref{eq::secular-equation} and where $\lambda_1$ is the smallest eigenvalue of $AA^T$, for all Newton iterates remain in \(\mathcal{L}\) and converge to \(\alpha^*\).
If \(\alpha = 0\) does not yield the solution, we initialize Newton's method from an \(\alpha_0\) chosen so that \(\phi'(\alpha_0)\) can be computed.
If ever $\alpha_k > \alpha^*$, convergence is not assured.
However, \citep[Lemma 7.3.3]{conn-gould-toint-2000} indicates that, due to the concavity of \(\phi\), the next iterate will either be in $\mathcal{L}$ or non-positive.
\Cref{eq:newton-safeguard} safeguards against \(\alpha < 0\) by restarting the Newton iterations from a positive value smaller than $\alpha_k$ in hopes of eventually landing in $\mathcal{L}$.

We now turn to the solution of~\eqref{eq:prox-problem-quadratic}.
For any $\alpha \geq 0$, the element \(u \in \prox{\nu \eta_\tau}(w, B)\) given by \Cref{cor:prox-quadratic} can equivalently be obtained from
\begin{equation}%
  \label{eq:augmented-matrix-pb}
  \begin{bmatrix}
    -Q & A^T      \\
    A  & \alpha I
  \end{bmatrix}
  \begin{bmatrix}
    u \\
    y
  \end{bmatrix}
  =
  \begin{bmatrix}
    -w \\
    -b
  \end{bmatrix}
  \quad \text{where} \ Q = \nu B + I
  .
\end{equation}
Because \(B\) will typically be a limited-memory quasi-Newton approximation in our implementation, materializing, let alone factorizing,~\eqref{eq:augmented-matrix-pb} would be inefficient.
Thus, we use MINRES-QLP \citep{choi-paige-saunders-2011}.
Even though only an approximate minimizer of~\eqref{eq:prox-problem-quadratic} is required, we compute an accurate solution of~\eqref{eq:augmented-matrix-pb} by setting both the absolute and relative tolerances of MINRES-QLP to \(\epsilon_M^{0.7}\) where \(\epsilon_M\) is the machine epsilon.
If \(\alpha = 0\) and \(A\) is rank deficient, MINRES-QLP is guaranteed to return the mininimum-norm solution of~\eqref{eq:augmented-matrix-pb}, but not necessarily \(y = -{(A Q^{-1} A^T)}^{\dagger} v\) required by \Cref{cor:prox-quadratic}.
Computing such \(y\) with MINRES-QLP can be done by solving a larger symmetric saddle-point system that represents the optimality conditions of the least-norm problem that defines \(y\).
Our preliminary implementation simply searches for \(\alpha \geq \sqrt{\varepsilon_M}\), where \(\varepsilon_M\) is the machine epsilon.
We will study and report on alternative computations of \(y\) in follow-up research.

In order to solve \(g(\alpha) = 0\) with $g$ as in \Cref{cor:prox-quadratic}, we solve the \emph{secular} equation
\begin{equation*}%
  \phi_Q(\alpha) = 0,
  \qquad
  \phi_Q(\alpha) \coloneqq \frac{1}{\|{(AQ^{-1} A^T + \alpha I)}^{-1} (AQ^{-1}x + b)\|_2} - \frac{1}{\nu\tau}.
\end{equation*}
If $y(\alpha)$ solves~\eqref{eq:augmented-matrix-pb},
\begin{equation}%
  \label{eq:newton-ratio}
  \frac{\phi_Q(\alpha)}{\phi_Q'(\alpha)} = \frac{\|y(\alpha)\|^2}{y{(\alpha)}^T{(AQ^{-1}A^T + \alpha I)}^{-1}y(\alpha)} \, (1 - \frac{\|y(\alpha)\|}{\nu\tau}).
\end{equation}
Computing the denominator above requires $w \coloneqq {(AQ^{-1}A^T + \alpha I)}^{-1}y(\alpha)$ found with MINRES-QLP by solving
\begin{equation}%
  \label{eq:augmented-matrix-pb-2}
  \begin{bmatrix}
    -Q & A^T      \\
    A  & \alpha I
  \end{bmatrix}
  \begin{bmatrix}
    u \\
    w
  \end{bmatrix}
  =
  \begin{bmatrix}
    0 \\
    y(\alpha)
  \end{bmatrix}
  .
\end{equation}









\section{Implementation and experiments}%
\label{sec:numerical}

We implemented all algorithms above in Julia $1.11$ \citep{bezanson-edelman-karpinski-shah-2017}.
\Cref{alg:exactpen,alg:R2N} are part of the \href{https://github.com/JuliaSmoothOptimizers/RegularizedOptimization.jl}{RegularizedOptimization.jl} package \citep{regularizedoptimization-jl}; \Cref{alg:proximal} is in \href{https://github.com/JuliaSmoothOptimizers/ShiftedProximalOperators.jl}{ShiftedProximalOperators.jl} \citep{shiftedproximaloperators-jl}.

We compare \Cref{alg:exactpen} against the augmented-Lagrangian solver Percival \citep{arreckx-lambe-martins-orban-2016,percival-jl} and the IPOPT \citep{wächter-biegler-2006,ipopt-jl} solver---which reduces to an SQP method for~\eqref{eq:nlp}---on equality-constrained problems from the CUTEst collection \citep{gould-orban-toint-2015,cutest-jl} with fewer constraints than variables, and with at most $300$ variables.
This results in a set of $50$ problems.
Each problem is given a limit of $5$ minutes of CPU time.
In what follows, we do not provide CPU time comparisons because our implementation is still preliminary, and memory usage is not optimal.

Percival and IPOPT are second-order methods while \Cref{alg:exactpen} with R2 is a first-order method.
For the comparison to be fair, in the problems that we expose to Percival and IPOPT, we replace the Hessian of the objective with a multiple of the identity, $\nabla^2 f(x) \approx \sigma_f(x) I$, and replace the constraint Hessians with the zero matrix, $\nabla^2 c_i(x) \approx 0$ for $1 \leq i \leq m$.
Because R2 updates \(\sigma_{k, j}\) at each iteration, we update \(\sigma_f(x)\) as in the spectral gradient method \citep{birgin-martinez-raydan-2014}, i.e., \(\sigma_f(x_{k, j}) \coloneqq s_{k, j}^T y_{k, j} / s_{k, j}^T s_{k, j}\),
where $s_{k, j} \coloneqq x_{k, j} - x_{k, j-1}$ and $y_{k, j} \coloneqq \nabla f(x_{k, j}) - \nabla f(x_{k, j-1})$.

The optimality measures in~\eqref{eq:xi} and~\eqref{eq:theta}, though appropriate in our context, are inconvenient to compare with solvers for smooth optimization, as the latter use criteria based on the KKT conditions for~\eqref{eq:nlp}.
Therefore, our implementation of \Cref{alg:exactpen} stops as soon as
\begin{equation}%
  \label{eq:stopping-criterion-ipopt}
  \|\nabla f(x_k) + J{(x_k)}^T y_k \|_\infty \leq \epsilon
  \quad \text{and} \quad
  \|c(x_k)\|_\infty \leq \epsilon,
\end{equation}
where \(y_k\) is the solution of~\eqref{eq::thProx-dual0} in the case where \(Q = \sigma_k I\).

\begin{figure}[ht]%
  \resizebox{.49\textwidth}{!}{%
    \includegraphics{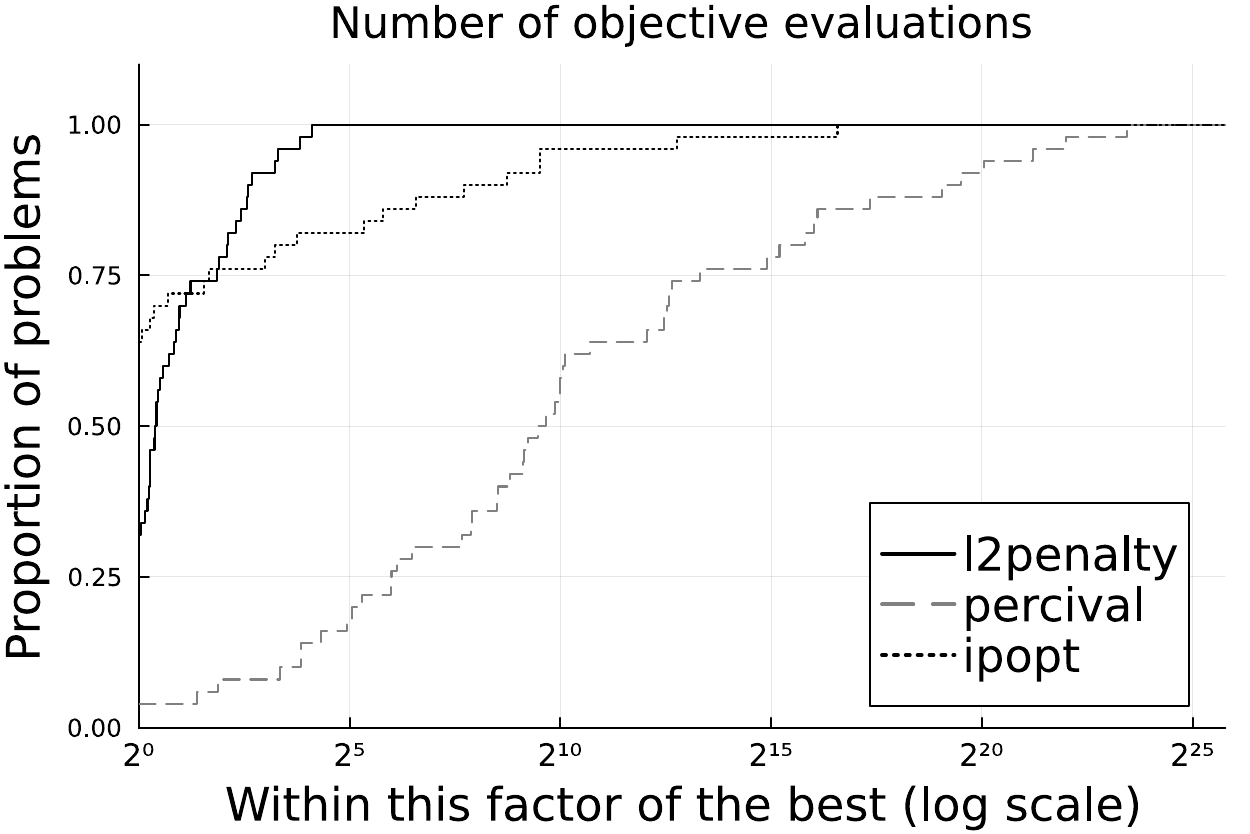}
  }
  \hfill
  \resizebox{.49\textwidth}{!}{%
    \includegraphics{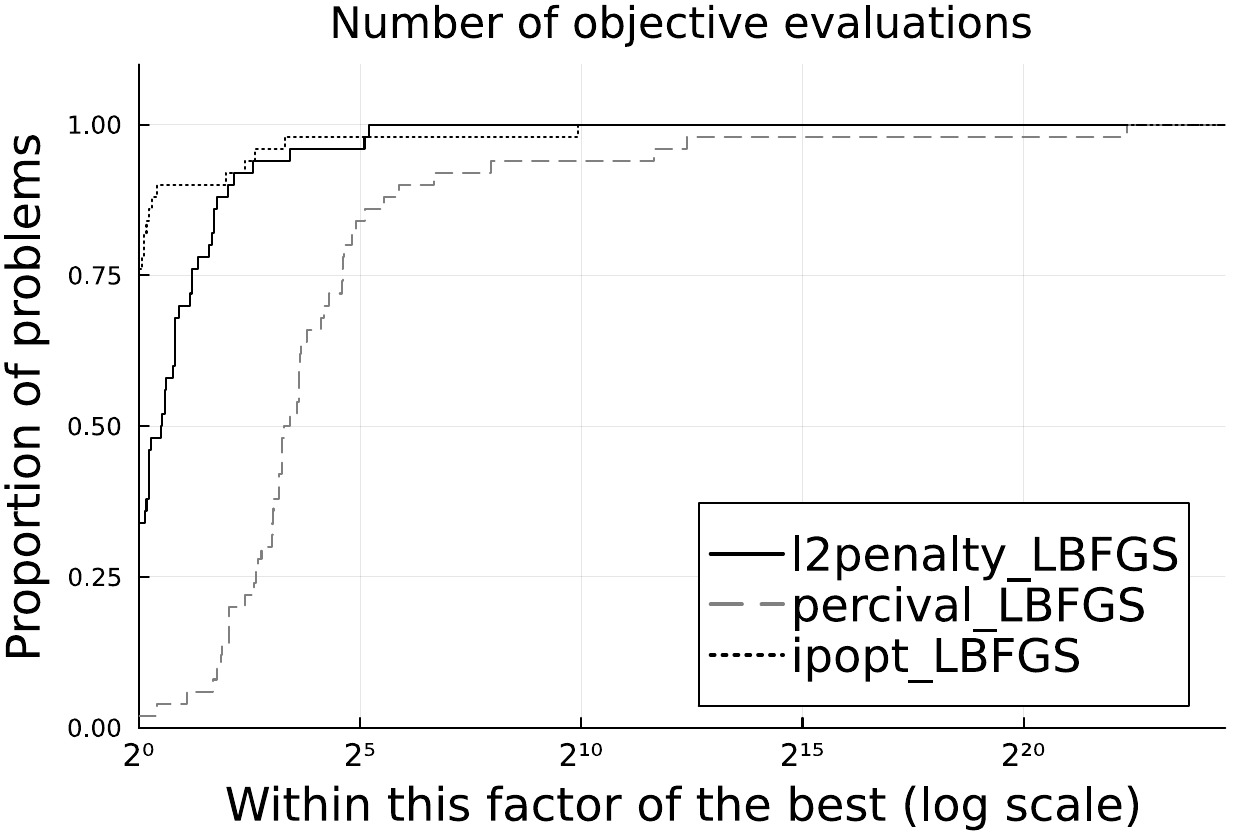}
  }\\
  \vfill
  \resizebox{.49\textwidth}{!}{%
    \includegraphics{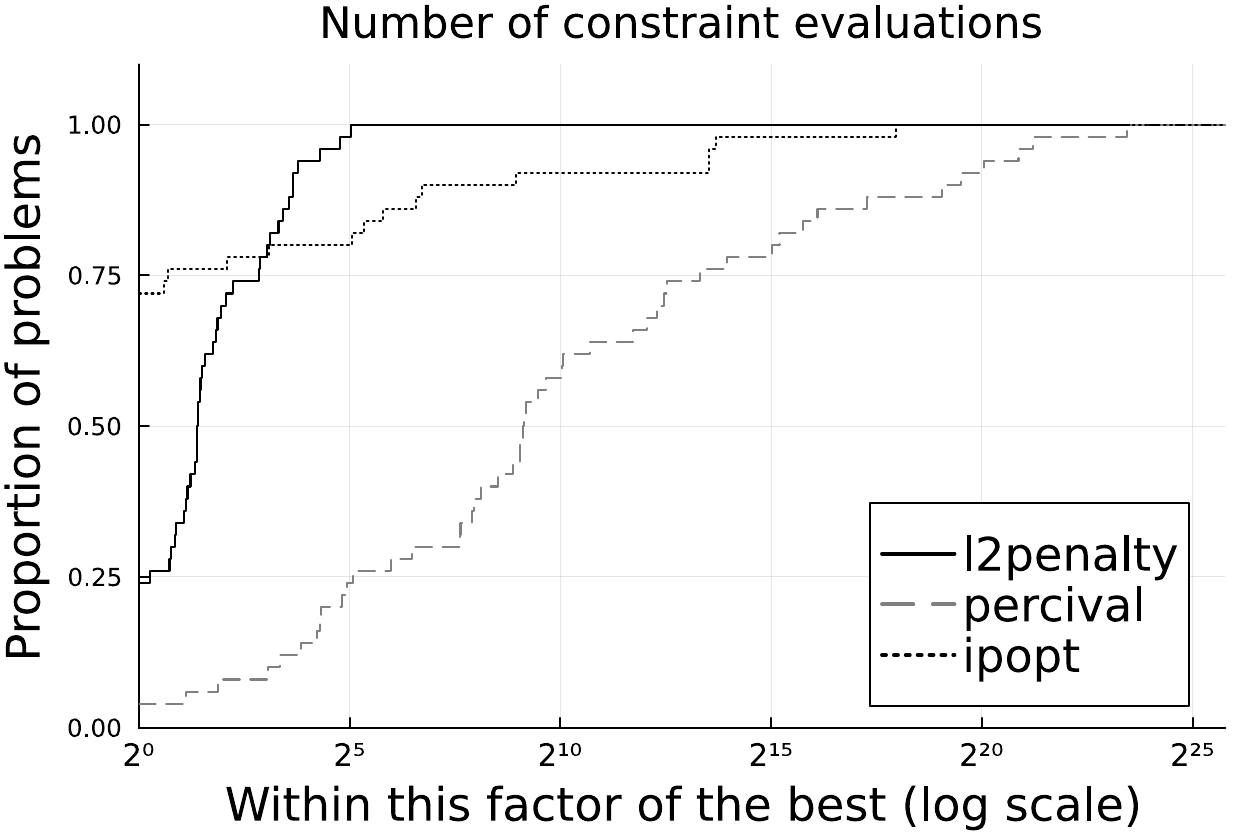}
  }
  \hfill
  \resizebox{.49\textwidth}{!}{%
    \includegraphics{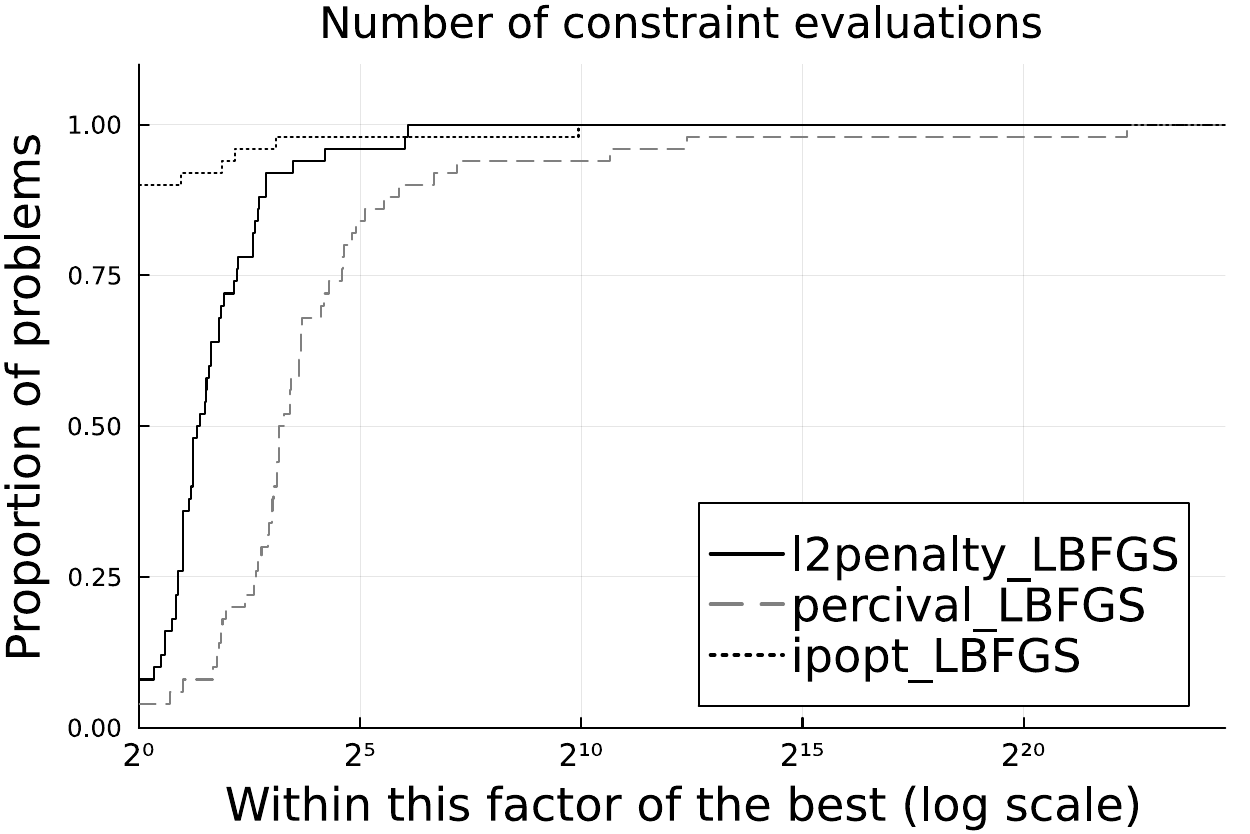}
  }\\
  \vfill
  \resizebox{.49\textwidth}{!}{%
    \includegraphics{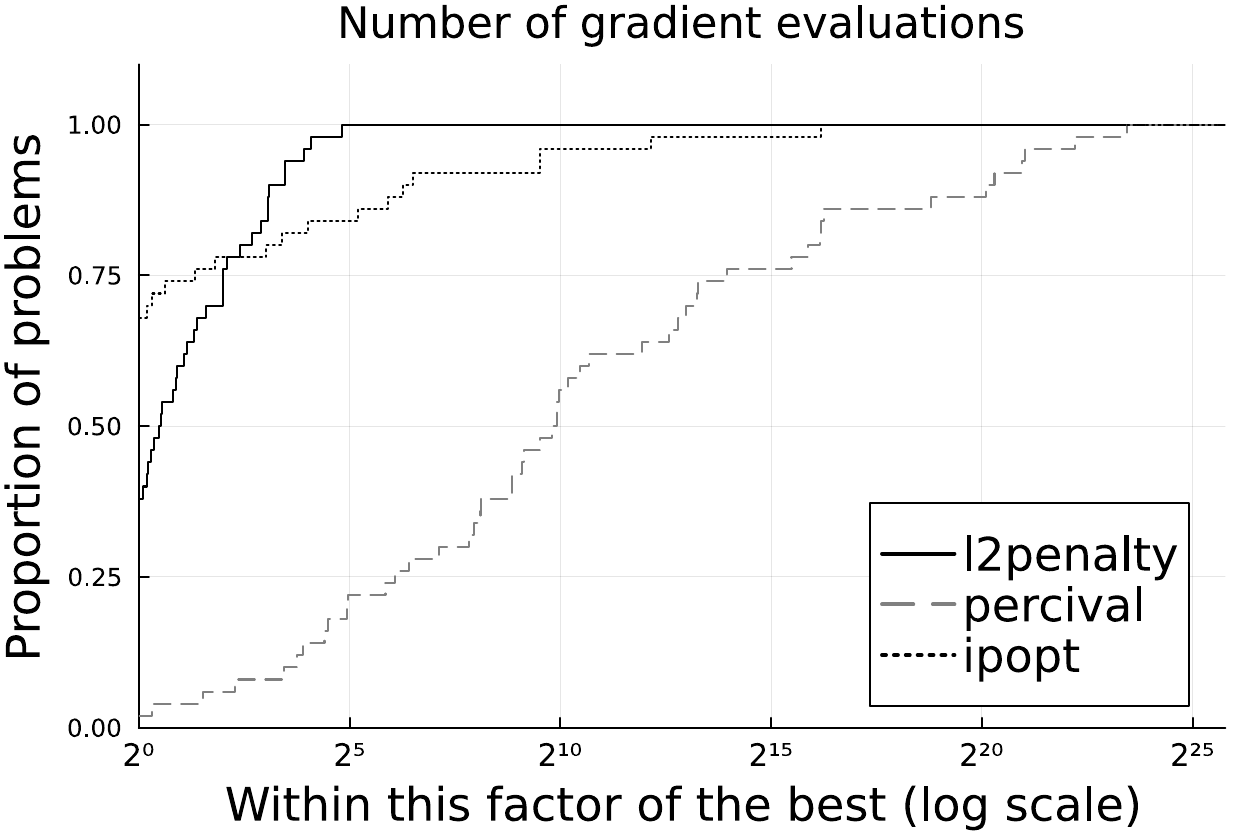}
  }
  \hfill
  \resizebox{.49\textwidth}{!}{%
    \includegraphics{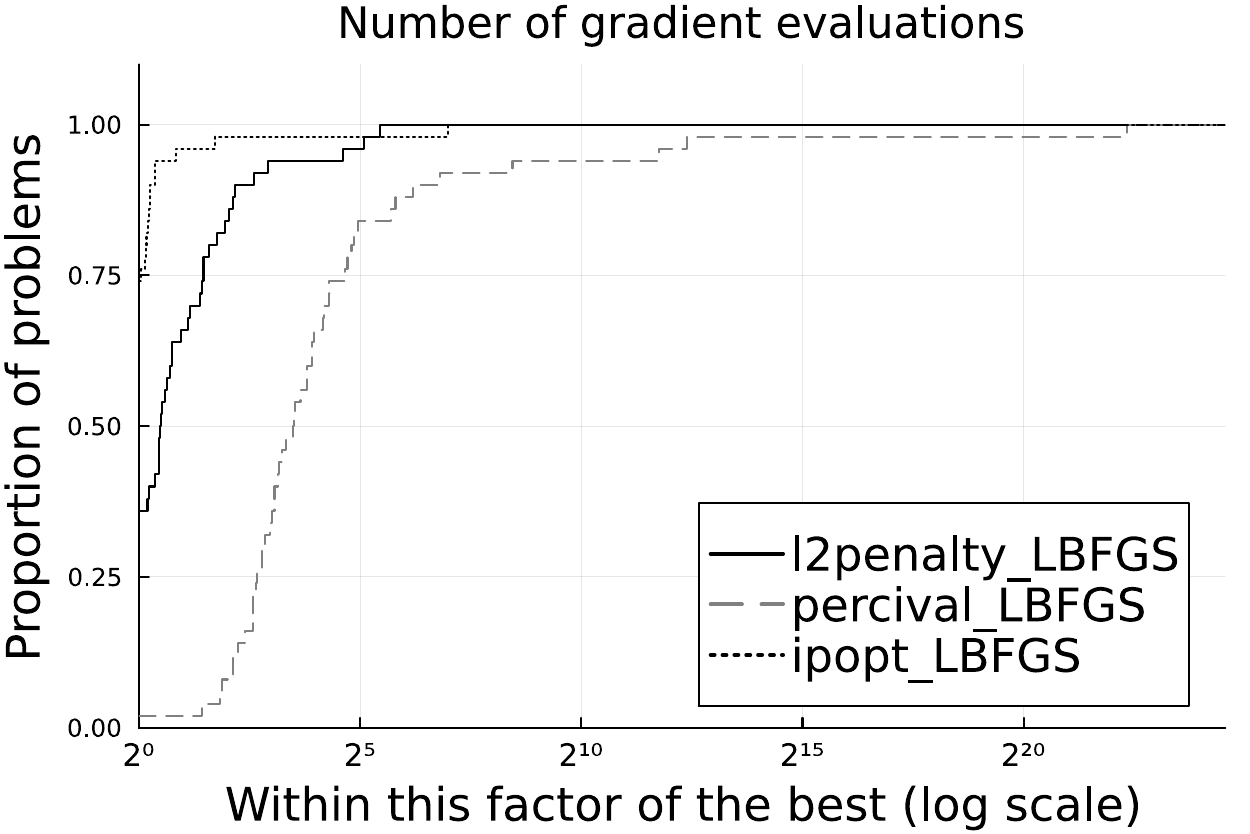}
  }
  \caption{Left: \Cref{alg:exactpen} with the R2 subsolver against Percival and IPOPT with a spectral quasi-Newton approximation.
    Right: \Cref{alg:exactpen} with the R2N subsolver against Percival and IPOPT with a LBFGS quasi-Newton approximation.}%
  \label{fig:perf-prof}
\end{figure}

Percival is modified so that is also stops when it finds points that satisfy~\eqref{eq:stopping-criterion-ipopt}, which is the criterion used in IPOPT\@.
We disabled IPOPT's automatic scaling of the problem statement, by setting \textit{nlp\_scaling\_method} to ``none’’, and the scaling of the dual feasibility with respect to the dual variables by setting \(s_{\max}\) to the maximum of the Float64 number system.
In these preliminary small-scale tests, all solvers use an absolute stopping test with $\epsilon = 10^{-3}$.
This is done in IPOPT by setting the parameter \textit{acceptable\_iter} to \(0\).
In \Cref{alg:exactpen}, we used $\beta_1 = \tau_0 = \sqrt{n m}$, $\epsilon_0 = 10^{-2}$, $\beta_2 = 0.1$, $\beta_3 = 10^{-2}$, $\beta_4 = \varepsilon_M \approx 2.2 \times 10^{-16}$, where \(n\) and \(m\) are respectively the number of variables and constraints of the problem.
In \Cref{alg:proximal}, we used \href{https://github.com/JuliaSmoothOptimizers/QRMumps.jl}{QRMumps.jl} \citep{qrmumps-jl} for the QR factorizations and set $\theta = 0.8$.
The stopping criterion for \Cref{alg:proximal} is $| \|s(\alpha_k)\|_2 - \nu\tau | < \varepsilon_M^{0.3}$.
We approximate the least-norm solution of $AA^T x = b$ using the Golub-Riley iteration \citep{dax-eldén-1998} with a regularization parameter of \(\varepsilon_M^{0.9}\).
The latter requires a single QR factorization and may be viewed as a form of iterative refinement.
More sophisticated strategies will be required in later versions of our implementation.
We set the maximal number of iterations for \Cref{alg:proximal} to $10$.
All tests are performed in double precision.

We also compare \Cref{alg:exactpen} combined with R2N to IPOPT and Percival.
All parameters are the same as before, and we experiment with both LBFGS and LSR1 Hessian approximations with memory \(6\).
At each successful iteration, the difference of gradients of the Lagrangian of~\eqref{eq:nlp},
\begin{equation*}%
  \nabla \mathcal{L}(x_k, y^{LS}_k) - \nabla \mathcal{L}(x_{k-1}, y^{LS}_k),
\end{equation*}
is used to update \(B_k\), where \(y^{LS}_k\) is the minimum least-squares solution of \(\nabla f(x_k) + {J(x_k)}^T y = 0\).
When we use an LBFGS approximation, our implementation uses damping~\citep[Section 5]{powell-1978} with damping parameter \(0.1\).
The MINRES-QLP tolerance is set to \(\varepsilon_M\) and its maximum number of iterations to $10^4$ to evaluate~\eqref{def:prox-quadratic}.

We report our results in the form of \citet{dolan-more-2002} performance profiles in \Cref{fig:perf-prof} in terms of number of objective, gradient and constraint evaluations.
Recall that a performance profile \citep{dolan-more-2002} is defined as follows.
For a set \(P\) of problems and a set \(S\) of solvers, we compute a performance metric \(t_{p, s}\) for each \(s \in S\) and \(p \in P\).
We then compute the ratio \(r_{p, s} \coloneqq t_{p, s}/\min_s t_{p,s} \) for each \(p \in P\) and \(s \in S\).
If \(n_p\) is the number of problems in \(P\), the performance profile for solver \(s\) is
\begin{equation*}
  \pi_s(\omega) \coloneqq \tfrac{1}{n_p} \sum_{p \in P} \mathbf{1}_{r_{p, s} \leq \omega},
\end{equation*}
where, for \(\omega \geq 0 \), \(\mathbf{1}_{r_{p, s} \leq \omega} = 1\) if \( r_{p, s} \leq \omega\), and \(\mathbf{1}_{r_{p, s} \leq \omega} = 0\) otherwise.

When using R2N, we only report results with LBFGS\@; results with LSR1 are nearly identical.
In the left column of \Cref{fig:perf-prof}, we see that our implementation is competitive with first-order Percival and IPOPT on all three metrics.
Because Percival never evaluates the Jacobian but only Jacobian-vector products, we do not provide profiles in terms of Jacobian evaluations.
The profiles suggest that the performance of \Cref{alg:exactpen} is very promising for problems where, for some reason, constraint Hessians cannot be evaluated or approximated.
In the right column of \Cref{fig:perf-prof}, we see that \Cref{alg:exactpen} is ahead of Percival, while IPOPT is competitive with our implementation in terms of objective and gradient evaluations.
Whether we use R2 or R2N, the robustness of \Cref{alg:exactpen} is on par with the other two solvers.
From the \(50\) problems considered in our experiment, two of them are infeasible, namely SSINE and VANDANIUMS\@.
Both where successfully detected as infeasible by our algorithm when using R2, whereas IPOPT only found VANDANIUMS and Percival found none.
When using R2N, only VANDANIUMS was detected by both IPOPT and \Cref{alg:exactpen}, and IPOPT reported a false positive with DIXCHLNG, Percival did not find any in this case either.
We think that these results are strongly in favor of considering exact penalty approaches as an alternative to augmented-Lagrangian approaches.

\section{Discussion and future work}%
\label{sec:discussion}

\Cref{alg:exactpen} is, to the best of our knowledge, the first practical implementation of the exact \(\ell_2\)-penalty method.
On small-scale problems, our preliminary implementation is ahead of an augmented Lagrangian method.
It is also ahead of IPOPT when spectral Hessian approximations are used in the latter in terms of efficiency and robustness.
It remains competitive in terms of robustness when using limited-memory quasi-Newton Hessian approximations.
Improvements in terms of efficiency are the subject of active research.

Among possible improvements, the solution of the secular equations requires attention in the presence of rank-deficient Jacobians when the solution is at, or near, the origin.
A possible improvement in terms of CPU time is to solve~\eqref{eq::secular-equation} in reduced-precision arithmetic.
Similarly, MINRES-QLP does not provide the required least-norm solution in \Cref{cor:prox-quadratic}.
When it appears that the solution is \(\alpha^* = 0\) and the Jacobian is rank deficient, a different saddle-point system can be solved that provides the required solution.

Our framework is general and allows us to implement the exact penalty method in any other norm as long as the resulting trust-region subproblem can be solved efficiently.
In particular, \Cref{th::thProx-dual} suggests that we could implement the exact $\ell_1$-penalty method originally proposed by \citet{pietrzykowski-1969}.
Indeed, evaluating the proximal operator amounts to solving a convex bound-constrained subproblem for which there are efficient polynomial time algorithms \citep[Chapter 7.8]{conn-gould-toint-2000}.

A trust-region variant of \Cref{alg:exactpen} based on \citep[Algorithm~\(3.1\)]{aravkin-baraldi-orban-2022} instead of R2 would enjoy the same asymptotic complexity bounds, although evaluation of the proximal operators would be significantly harder.

In future research, we plan to extend our analysis to constrained problems where the objective is the sum of a smooth function $f$ and a nonsmooth regularizer $h$:
\begin{equation*}
  \minimize{x} f(x) + h(x) \ \st \ c(x) = 0.
\end{equation*}

\small
\subsection*{Acknowledgements}
We would like to thank Tangi Migot and Alexis Montoison for their kind assistance with the numerical implementation of our algorithm.
\normalsize

\appendix
\section{Inner algorithms}%
\label{sec:inner_algs}

\Cref{alg:R2N} states \citep[Algorithm~\(4.1\)]{diouane-habiboullah-orban-2024} for completeness; \citep[Algorithm~\(6.1\)]{aravkin-baraldi-orban-2022} is a special case of \Cref{alg:R2N} with \(B_j = 0\) and \(\alpha = 1\).

\begin{algorithm}[ht]%
  \caption{%
    \label{alg:R2N}
    R2N\@: A proximal modified Quasi-Newton method, \citep[Algorithm 4.1]{diouane-habiboullah-orban-2024}
  }
  \begin{algorithmic}[1]%

    \Require{}%
    The value \(\tau \coloneqq \tau_k > 0\) of the penalty parameter is given by the outer iteration.

    \State{}%
    Choose constants \(\sigma_{\min} > 0\), $0 < \eta_1 \leq \eta_2 < 1$ and $ 0 < \gamma_3 \leq 1 < \gamma_1 \leq \gamma_2$, \(0 < \alpha < 1\) and $\beta > 1$.%

    \State{}%
    Choose $x_0 \in \R^n$ where $h$ is finite, $\sigma_0 \geq \sigma_{\min}$, compute $f(x_0) + \tau h(x_0)$.%

    \For{$j = 0, 1, \ldots$}%
    \State{}%
    Choose $B_j = B_j^T \in \R^{ n \times n }$.%

    \State{}%
    Set $\nu_j = \alpha / (\|B_j\|_2 + \sigma_j)$.%

    \State{}%
    Compute $s_{j,\mathrm{cp}} \in \prox{\nu_j \tau\psi(\cdot; x_{j})}(-\nu_j\nabla f(x_j))$.%

    \State{}%
    Compute an approximate minimizer $s_j$ of \(m_Q(s; x_j, B_j) + \tfrac{1}{2} \sigma_j \|s\|^2_2\) such that
    \begin{equation*}
      m_Q(s_j; x_j, B_j) + \tfrac{1}{2} \sigma_j \|s_j\|^2_2 \leq m_Q(s_{j, \mathrm{cp}}; x_j, B_j) + \tfrac{1}{2} \sigma_j \|s_{j, \mathrm{cp}}\|^2_2.
    \end{equation*}

    \State{}%
    \label{eq:R2N-unbounded-check}%
    If $\| s_j \|_2 > \beta \| s_{j,\mathrm{cp}} \|_2$, reset $s_j = s_{j,\mathrm{cp}}$.%

    \State{}%
    Compute the ratio
    \begin{equation*}%
      \rho_j \coloneqq \frac{f(x_j) + \tau h(x_j) - \left( f(x_j + s_j) + \tau h(x_j + s_j) \right) }{\varphi_Q(0 ; x_j, B_j) + \tau\psi(0 ; x_j) - \left( \varphi_Q(s_j ; x_j, B_j) + \tau\psi(s_j ; x_j) \right)}.
    \end{equation*}
    \State{}%
    If $\rho_j \geq \eta_1$, set $x_{j+1} = x_j + s_j$.
    Otherwise, set $x_{j+1} = x_j$.
    \State{}%
    Update the regularization parameter according to
    \begin{equation*}
      \sigma_{j+1} \in
      \begin{cases}
        \left[ \gamma_3\sigma_j, \sigma_j \right]         & \text{if} \ \rho_j \geq \eta_2,          \\
        \left[ \sigma_j, \gamma_1\sigma_j \right]         & \text{if} \ \eta_1 \leq \rho_j < \eta_2, \\
        \left[ \gamma_1\sigma_j, \gamma_2\sigma_j \right] & \text{if} \ \rho_j < \eta_1.
      \end{cases}
    \end{equation*}
    \State{}%
    Set \(\sigma_{j+1} \leftarrow \max(\sigma_{j+1}, \, \sigma_{\min})\).%
    \EndFor%
  \end{algorithmic}
\end{algorithm}

Since $h$ is bounded below, finding an $x_0$ where $h$ is finite in \Cref{alg:R2N} is always trivial.
Although neither \citet{aravkin-baraldi-orban-2022} nor \citet{diouane-habiboullah-orban-2024} mention a minimum regularization parameter \(\sigma_{\min}\), it does not change the complexity of their algorithms.
For example, building upon \citep{aravkin-baraldi-orban-2022}, \citet{leconte-orban-2025} use a maximum trust-region parameter \(\Delta_{\max}\) that is analogous to \(\sigma_{\min}\) in our context and which does not change the complexity of the algorithm.
In view of \Cref{th::thProx-dual,th::thProx-dual-p2}, \(B_j + \sigma_j I\) should be positive definite when computing a step.
We argue however that if this is not the case then \Cref{eq:R2N-unbounded-check} will be triggered, causing the algorithm to retreat to the Cauchy point

\begin{lemma}%
  \label{lemma:positive-def}
  Let $d \in \R^n$, $B = B^T \in \R^{n \times n}$, $A \in \R^{m \times n}$, $b \in \R^m$, $\tau > 0$, and $\sigma > 0$.
  If $B + \sigma I$ is not positive semi-definite,
  \begin{equation}%
    m(s) \coloneqq d^T s + \tfrac{1}{2} s^T B s + \tau \|As + b\|_2 + \tfrac{1}{2} \sigma \|s\|^2_2
  \end{equation}
  is such that \(\lim_{t \to \infty} \min_{\|s\| \leq t} m(s) = -\infty\).
\end{lemma}

\begin{proof}
  Assume that $B + \sigma I$ is not positive semi-definite.
  Then, there is a vector $\bar{s} \in \R^n$ with $\| \bar{s} \|_2 = 1$ such that $\bar{s}^T (B + \sigma I) \bar{s} < 0$.
  This implies that
  \begin{align*}
    \lim_{t \to \infty} m(t \bar{s}) \leq \lim_{t \to \infty} \tfrac{1}{2} t^2 \bar{s}^T (B + \sigma I) \bar{s} + t (d^T \bar{s}^T + \tau \|A\bar{s}\|_2) + \tau \|b\|_2 = -\infty.
  \end{align*}
\end{proof}

From \Cref{lemma:positive-def}, whenever \(B_j + \sigma_j I\) is not positive definite, the minimizer of \(m_Q(\cdot; x_j, B_j) + \tfrac{1}{2}\sigma_j \|\cdot\|^2_2\) is arbitrarily large so that \Cref{eq:R2N-unbounded-check} is activated for any \(\beta > 0\).



\bibliographystyle{abbrvnat}
\bibliography{abbrv,l2-penalty-siam}


\end{document}